\documentclass{amsart}

\usepackage{amssymb}
\usepackage{mathtools}
\usepackage{bm}
\usepackage{enumitem}
\usepackage{microtype}
\usepackage[hidelinks]{hyperref}
\hypersetup{
  pdftitle={Weighted Harris--oscillation contraction and invariant expectations for nonlinear Markov semigroups},
  pdfauthor={Danyang Zhao and Huaizhong Zhao},
  pdfsubject={Invariant expectations and ergodicity for nonlinear Markov semigroups and G-SDEs},
  pdfkeywords={sublinear expectation, G-SDE, Harnack inequality, weighted oscillation, invariant expectation, ergodicity}
}

\numberwithin{equation}{section}

\theoremstyle{plain}
\newtheorem{theorem}{Theorem}[section]
\newtheorem{proposition}[theorem]{Proposition}
\newtheorem{lemma}[theorem]{Lemma}
\newtheorem{corollary}[theorem]{Corollary}
\newtheorem{assumption}{Assumption}
\theoremstyle{definition}
\newtheorem{definition}[theorem]{Definition}

\theoremstyle{remark}
\newtheorem{remark}[theorem]{Remark}

\newcommand{\R}{\mathbb R}

\newcommand{\Et}{\widetilde{\mathbb E}}
\newcommand{\Borel}{\mathcal B}
\newcommand{\Q}{\mathcal Q}
\newcommand{\osc}{\operatorname{osc}}
\newcommand{\tr}{\operatorname{tr}}

\newcommand{\cE}{\mathcal E}

\newcommand{\cL}{\mathcal L}

\newcommand{\norm}[1]{\left\lVert #1\right\rVert}
\newcommand{\abs}[1]{\left\lvert #1\right\rvert}
\newcommand{\Prog}{\operatorname{Prog}}

\title[Weighted Harris--oscillation contraction]{Weighted Harris--oscillation
contraction and invariant expectations for nonlinear Markov semigroups}

\author[D. Zhao
]{Danyang Zhao
}
\address{School of Mathematics, Shandong University, Jinan 250100, China
}
\email{dyang.z@mail.sdu.edu.cn
}

\author[H. Zhao
]{Huaizhong Zhao
}
\address{Department of Mathematical Sciences, Durham University,
Durham DH1 3LE, United Kingdom; and School of Mathematics, Shandong University, Jinan 250100, China
}
\email{huaizhong.zhao@durham.ac.uk
}

\date{}



\begin{document}

\begin{abstract}
We establish a weighted contraction principle for sublinear Markov
semigroups on noncompact state spaces.  Its abstract inputs are a Lyapunov
drift condition and a strict local contraction of oscillations on a
Lyapunov sublevel set.  For locally uniformly parabolic diffusion
semigroups, the latter is obtained by applying the parabolic Harnack
inequality to both the upper evolution and its conjugate lower evolution:
if $u=T_t\varphi$, then $u-m$ and $M-u$ solve the upper and lower Bellman
equations, respectively.  Combining the resulting local oscillation
contraction with the Lyapunov drift yields a strict contraction in a
weighted oscillation seminorm, a unique invariant sublinear expectation in
the finite-Lyapunov-moment class, and weighted exponential convergence for
Borel observables controlled by the Lyapunov function.

To treat bounded Borel data for $G$-SDEs, we use a canonical
covariance-law representation, progressively measurable state maps
compatible with conditioning and pasting, and a stopped dynamic programming
principle.  A probabilistic Krylov--Safonov argument then gives
positive-time H\"older regularisation before viscosity identification,
thereby avoiding a measurable-data envelope argument.  The abstract result
is verified for ordinary locally uniformly elliptic SDEs, multiplicative
$G$-SDEs, and the full equation with quadratic-variation drift
\[
 dX_t=b(X_t)\,dt
      +\sum_{i,j=1}^m h_{ij}(X_t)\,d\langle B^i,B^j\rangle_t
      +\sum_{j=1}^m\sigma_j(X_t)\,dB_t^j.
\]
For the full equation the hypotheses are expressed through the effective
drift and covariance and are required uniformly over all covariance
scenarios. Finally, the
contraction and stochastic continuity of the associated stationary process
yield ergodicity in the sense of Feng and Zhao.
\end{abstract}

\maketitle

\section{Introduction}

Nonlinear expectations provide a probabilistic framework for situations in
which uncertainty cannot be described adequately by a single reference
probability measure.  A sublinear expectation may be viewed as an upper
envelope of linear expectations over a family of plausible laws and thus
captures a robust, worst-case assessment of random outcomes.  Within this
framework, $G$-expectation and $G$-Brownian motion give a dynamically
consistent model of volatility uncertainty: instead of prescribing a single
covariance matrix, they allow the covariance to vary within a specified set.
$G$-diffusions are consequently nonlinear analogues of classical It\^o
diffusions and arise naturally in stochastic analysis under model
uncertainty; see Denis et al.\ \cite{DenisHuPeng2011} and Peng
\cite{Peng2004,Peng2005,Peng2007,Peng2019}.  For the classical
stochastic-control and viscosity-solution background, see Fleming and Soner
\cite{FlemingSoner2006}.  The backward theory, including nonlinear
Feynman--Kac representations and BSDEs driven by $G$-Brownian motion, was
developed by Hu et al.\ \cite{HuJiPengSong2014}.  The same
volatility-uncertainty paradigm underlies the theory of ambiguous volatility
in mathematical finance Epstein and Ji \cite{EpsteinJi2013}.

The long-time dynamics of these systems are important because they determine
whether the effect of the initial state is eventually lost and whether the
uncertain system possesses a well-defined stationary regime.  An invariant
sublinear expectation plays the role of an invariant probability measure,
while convergence toward it quantifies robust stability and mixing.
Ergodicity further rules out a nontrivial decomposition into invariant
components.  In the sublinear expectation formulation used below, an invariant sublinear
expectaion $\mathbb E$ is in the $\left(0,0\right)$-regime when every
invariant event $A$ satisfies
\[
 \mathbb E(I_A)=0
 \quad\text{or}\quad
 \mathbb E(I_{A^c})=0;
\]
see Feng and Zhao \cite{FengZhao2021}.  These properties form the basis for
stationary analysis and long-run asymptotic questions under volatility
uncertainty.  They are genuinely nonlinear issues, because there need not be
a single transition kernel or a single invariant law to which classical
Markov-process theory can be applied directly; see also Hu et al.
\cite{HuLiWangZheng2015}, and Zhao and
Zhao \cite{ZhaoZhao2026}.

There is also a closely related ergodic theory for invariant capacities.
Cerreia-Vioglio et al.\ \cite{ergodiclower} established an
ergodic theorem for continuous lower probabilities and studied the
$\left(0,1\right)$-regime, in which every invariant set has capacity either
zero or one.  Sheng and Song \cite{ergodicuppersong} obtained a structural
characterisation of continuous ergodic upper probabilities and strengthened
the corresponding ergodic theorem.  More recently, Feng et al.
developed Birkhoff, subadditive and multiplicative ergodic theorems, together
with a theory of mixing, for invariant capacities
\cite{fenghuangliuzhao2026}.  For invariant upper probabilities, the
$\left(0,1\right)$-regime does not in general imply irreducibility.  Feng et
al.\ \cite{fenghuangliuzhao2024} proved that it is equivalent
to their finite-ergodic-components property.  This contrasts with the
$\left(0,0\right)$-regime, which is the
irreducibility-type notion used in the present paper.

Related work on long-time behaviour under probability uncertainty includes
the theory of imprecise Markov chains.  For finite-state chains, de Cooman et
al.\ \cite{imc2009} established conditions ensuring
convergence, independently of the initial uncertainty model, to a unique
invariant upper expectation.  T'Joens and De Bock \cite{imc2021}
subsequently introduced a notion of weak ergodicity based on convergence of
upper and lower expected time averages and obtained a necessary and
sufficient characterisation in terms of graph-theoretic accessibility
conditions.

The classical theory of Markov-process stability rests on two ingredients.
Under the usual irreducibility and petite-set hypotheses, Foster--Lyapunov
drift criteria control excursions and yield positive Harris recurrence and
the existence of an invariant measure Meyn and Tweedie 
\cite{MeynTweedie1993,MeynTweedie2009}.  A small-set minorisation supplies
the complementary local mixing mechanism.  Down et al.
\cite{DownMeynTweedie1995} characterise exponential ergodicity through an
appropriate combination of drift and minorisation conditions.  Abstracting
this picture gives the Harris principle: a Lyapunov condition drives the
dynamics back toward a bounded region, while a local mixing condition
contracts discrepancies between states once they enter that region.  The
weighted-norm formulation of Hairer and Mattingly
\cite{HairerMattingly2011} shows particularly clearly how the local and
global mechanisms combine to produce a strict contraction, a unique
invariant measure, and an exponential convergence estimate.  For sublinear
Markov semigroups, however, transition probabilities are replaced by upper
expectations over families of laws, and a common minorising measure may be
unnatural or difficult to verify.

Under the loss convention, nonlinear expectations also admit a natural
interpretation in terms of coherent risk measures.  If $\varphi$ represents
a terminal loss, then $T_t\varphi(x)$ is the worst-case expected loss over
the probabilistic scenarios encoded by the $G$-expectation; see Artzner et
al.\ \cite{ArtznerDelbaenEberHeath1999} and Delbaen
\cite{Delbaen2002}.  The semigroup property expresses the dynamic
consistency of these evaluations, in the spirit of dynamic coherent risk
measures; see Riedel \cite{Riedel2004}.

The long-time problem is important from this perspective.  A finite-horizon
risk evaluation generally depends on the present state of the system,
whereas long-horizon planning requires a canonical assessment that is stable
under the uncertain dynamics.  Our convergence theorem shows that the
dependence on the initial state is forgotten exponentially fast and that the
worst-case expected loss approaches the invariant sublinear expectation.
The latter therefore provides a stationary long-run risk benchmark under
volatility uncertainty.  The precise invariant-risk identity, together with
a possible entrance-family counterpart for nonautonomous evolutions, is
recorded in Remark~\ref{rem:long-run-risk}.

Invariant nonlinear expectations for $G$-diffusions have been studied under
strong dissipativity assumptions; see, for example, Hu et al.
\cite{HuLiWangZheng2015}.  Feng and Zhao \cite{FengZhao2021} developed the
corresponding ergodic theory of sublinear Markovian processes through a
canonical dynamical system on a lifted sublinear expectation space.  They
also treated $G$-Brownian motion on the circle, providing a compact-state
example outside the strongly dissipative Euclidean framework.  That example
did not, however, resolve the general problem for $G$-SDEs on $\mathbb R^d$
under weak dissipativity and ellipticity.  More recently, Zhao and Zhao
\cite{ZhaoZhao2026} developed a nonlinear Harris theorem based on a local
coupling method and applied it to additive $G$-SDEs, thereby providing a
partial answer.

The present paper uses a different local input.  We work with contraction of
oscillations and derive it from two parabolic Harnack inequalities, one for
the sublinear semigroup and one for its conjugate lower semigroup.  This
isolates the precise part of the Harris argument needed by the abstract
nonlinear theory and enlarges the class of diffusion models for which the
local hypothesis can be verified.  In particular, it allows us to treat
multiplicative-noise $G$-SDEs under weak dissipativity and local uniform
ellipticity.

The basic observation is simple but decisive in the sublinear setting.  Let
\[
 u(t,x)=T_t\varphi(x),
 \qquad
 m=\inf_{\mathfrak Q}u,
 \qquad
 M=\sup_{\mathfrak Q}u
\]
on a suitable space--time cylinder $\mathfrak Q$.  For a linear diffusion,
$u-m$ and $M-u$ solve the same equation.  For an upper nonlinear semigroup,
by contrast, $u-m$ solves the upper Bellman equation, whereas $M-u$ solves
the conjugate lower Bellman equation.  Harnack must therefore be applied to
both equations.  Adding the resulting estimates recovers a fixed fraction
of $M-m$ at a later time and yields the strict local oscillation contraction
in Lemma~\ref{lem:Hosc}.  A one-sided Harnack estimate alone does not produce
the full oscillation contraction.

The abstract theorem is consequently stated directly in terms of local
oscillation contraction rather than Harnack regularity.  If $V$ is a
Lyapunov function and $C_R=\{V\leq R\}$, its hypotheses combine a drift
estimate at some time $a>0$ with an inequality of the form
\[
 \osc_{C_R}(T_a\varphi)
 \leq \vartheta_R\,
       \osc_{\mathfrak Q_R}(T_\cdot\varphi),
 \qquad \vartheta_R<1.
\]
Theorem~\ref{thm:weighted} turns these assumptions into
\[
 [T_a\varphi]_{\beta,V}
 \leq q[\varphi]_{\beta,V},
 \qquad q<1,
\]
for a weighted oscillation seminorm on functions modulo constants.  The
dual formulation gives a contraction on a complete metric space of
sublinear expectations having finite $V$-moment.  It follows that there is
a unique invariant expectation in this class, together with pointwise,
operator, and weighted exponential convergence; see
Corollary~\ref{thm:invariant}.  Thus the nonlinear fixed-point argument is
fully separated from the analytic or probabilistic verification of the
local oscillation estimate.

There is a second difficulty specific to the $G$-SDE applications.  The
natural $G$-expectation space does not by itself assign values to all bounded
Borel path functionals, whereas invariant-expectation and Harris arguments
are most naturally formulated on Borel state functions.  We therefore
specify the extension through the canonical family of covariance laws and
use the conditioning-and-pasting construction of Nutz and van Handel
\cite{NutzVanHandel2013}, together with the quasi-sure aggregation viewpoint
of Soner et al.\ \cite{SonerTouziZhang2011}.  For the full $G$-SDE,
a universal progressively measurable version of the state process is
constructed through canonical Euler approximations.  This supplies the
canonical flow and stopped dynamic programming principle simultaneously
across the nondominated covariance laws, without selecting a separate
Wiener-space realisation for each law.

The order of the regularity argument is essential.  The Bellman statement
is first proved for bounded uniformly continuous initial data.  For bounded
Borel data, Lemma~\ref{lem:BorelDataSmoothing} uses the stopped dynamic
programming principle and the probabilistic Krylov--Safonov hitting estimate
\cite{Krylov1987,KrylovSafonov1979} to establish positive-time local
H\"older continuity directly.  Once this continuity is known, the usual
dynamic-programming test-function argument identifies the Borel-data value
function as a viscosity solution of the upper Bellman equation; applying
the analogous argument to a constant minus the value function gives the
conjugate lower equation.  This architecture, used in
Proposition~\ref{prop:multbellman} and
Corollary~\ref{cor:multregularisation}, and again in
Proposition~\ref{prop:fullHarnack} and
Corollary~\ref{cor:fullregularisation}, avoids a circular appeal from
Bellman regularity to Borel-data smoothing and avoids the need for a
separate measurable-data envelope theorem.

The main contributions may be summarised as follows.
\begin{enumerate}[label=\textup{(\roman*)},leftmargin=*]
 \item We prove the weighted nonlinear Harris--oscillation theorem,
 Theorem~\ref{thm:weighted}.  The abstract local hypothesis is precisely the
 oscillation estimate used in the proof, so the theorem also applies when
 that estimate is obtained by a method other than Harnack.

 \item We give a two-sided Harnack mechanism for locally uniformly
 parabolic upper and lower Bellman equations and combine it with a
 Borel-data smoothing argument based on stopped dynamic programming.  The
 resulting verification is uniform over the admissible covariance laws and
 does not presuppose a common transition-density minorisation.

 \item We apply the framework first to ordinary locally uniformly elliptic
 SDEs as a linear benchmark and then to multiplicative and full $G$-SDEs.
 For the full equation, each covariance scenario $Q\in\Q$ has effective
 coefficients
 \[
  b_Q(x)=b(x)+\sum_{i,j=1}^m Q^{ij}h_{ij}(x),
  \qquad
  a_Q(x)=\sigma(x)Q\sigma(x)^{\mathsf T}.
 \]
 The assumptions are a robust Lyapunov inequality and local uniform
 ellipticity for these coefficients, both uniform in $Q$.  Theorems
 \ref{thm:multG} and \ref{thm:fullG} then give unique invariant sublinear
 expectations and weighted exponential convergence.

 \item We connect the contraction result with the ergodicity criterion of
 Feng and Zhao \cite{FengZhao2021}.  Convergence makes the fixed space
 simple, while a direct moment estimate proves stochastic continuity of the
 canonical stationary process.  This yields Corollary~\ref{cor:ergodic}
 without imposing a mean-uncertainty condition for the forward implication.
\end{enumerate}

It is useful to emphasise what is, and is not, asserted in the $G$-SDE
verification.  Under a law in the canonical covariance family, the
quadratic-variation density may be progressively chosen.  The Bellman
operator records pointwise covariance scenarios through $b_Q$ and $a_Q$,
but the proof never replaces an adapted volatility process by a single
deterministic coefficient having the same transition operator.  Nor does it
seek a minorising probability measure common to all admissible diffusions.
The common structure is instead the pair of Pucci inequalities for the upper
and lower equations, from which the Harnack constants are uniform over the
covariance scenarios.

In the present paper, we restrict attention to globally Lipschitz drift
coefficients, which in particular have at most linear growth.  The abstract
weighted nonlinear Harris--oscillation contraction principle is not
intrinsically tied to this growth restriction.  We expect that, when
combined with suitable localisation, non-explosion, and moment estimates,
the method can be extended to locally Lipschitz, weakly dissipative drifts of
polynomial growth.  Carrying out this extension requires additional
technical arguments concerning the canonical state map, dynamic
programming, and continuity with respect to the initial state.  To keep the
present paper focused and of manageable length, we do not pursue this
extension here.

The paper is organised as follows.  Section~\ref{sec:prelim} introduces
sublinear Markov semigroups, weighted function spaces, and the metric on
sublinear expectations.  Section~\ref{sec:harnack} converts the two-sided
Harnack property into local oscillation contraction, and
Section~\ref{sec:abstract} proves the abstract weighted contraction and
invariant-expectation results.  Section~\ref{sec:classical} treats ordinary
uniformly elliptic SDEs as a linear benchmark.
Section~\ref{sec:BorelDataSmoothing} establishes a general Borel-data
regularisation lemma for uniformly elliptic families satisfying a stopped
dynamic programming principle; this lemma is used in
Sections~\ref{sec:multG} and~\ref{sec:fullG}.
Section~\ref{sec:multG} develops the canonical dynamic programming, Bellman,
regularisation, and invariant-expectation theory for multiplicative
$G$-SDEs.  Section~\ref{sec:fullG} carries out the corresponding analysis for
the full equation with quadratic-variation drift.  Finally,
Section~\ref{sec:ergodicity} proves ergodicity of the associated canonical
stationary systems.

\section{Sublinear semigroups and weighted spaces}
\label{sec:prelim}

Let $(E,\Borel(E))$ be a measurable space.  In applications $E=\R^d$.
Let $B_b(E)$ denote the bounded Borel functions on $E$.
Let us recall some definitions first. 

\begin{definition}
A functional $\cE$ on a linear lattice of real functions containing the
constants is a \emph{sublinear expectation} if, for all admissible
$\varphi,\psi$ and all $\lambda\geq0$,
\begin{enumerate}[label=\textup{(\roman*)}]
 \item $\varphi\leq\psi$ implies $\cE[\varphi]\leq\cE[\psi]$;
 \item $\cE[c]=c$ for every $c\in\R$;
 \item $\cE[\varphi+\psi]\leq\cE[\varphi]+\cE[\psi]$;
 \item $\cE[\lambda\varphi]=\lambda\cE[\varphi]$.
\end{enumerate}
The conjugate lower expectation is
\[
 \underline{\cE}[\varphi]:=-\cE[-\varphi].
\]
\end{definition}

Translation invariance follows from the axioms:
\begin{equation}\label{eq:translation}
 \cE[\varphi+c]=\cE[\varphi]+c.
\end{equation}
Moreover,
\begin{equation}\label{eq:absbound}
 -\cE[\abs{\varphi}]\leq \cE[\varphi]
 \leq \cE[\abs{\varphi}].
\end{equation}

\begin{definition}
A family $(T_t)_{t\geq0}$ of maps on $B_b(E)$ is a \emph{sublinear Markov
semigroup} if, for every $t\geq0$ and $x\in E$, the map
$\varphi\mapsto T_t\varphi(x)$ is a sublinear expectation and
\[
 T_0=I,\qquad T_{t+s}=T_tT_s,\qquad s,t\geq0.
\]
Its conjugate lower semigroup is
\[
 \underline T_t\varphi=-T_t(-\varphi).
\]
\end{definition}

For every bounded $\varphi$, define 
\[\osc_{E}(\varphi)=\sup_{x\in E}\varphi(x)-\inf_{x\in E}\varphi(x).\]
It is easy to see that
\begin{equation}\label{eq:osc-nonexp}
 \inf_E\varphi\leq \underline T_t\varphi
 \leq T_t\varphi\leq\sup_E\varphi,
 \qquad
 \osc_E(T_t\varphi)\leq\osc_E(\varphi).
\end{equation}

\begin{definition}
A sublinear expectation $\Et$ is \emph{invariant} for $(T_t)$ if
\begin{equation}\label{eq:invdef}
 \Et[T_t\varphi]=\Et[\varphi],\qquad t\geq0.
\end{equation}
\end{definition}

Let $V:E\to[0,\infty)$ be measurable and finite.  Define the weighted Borel function space
\[
 B_V(E)
 :=
 \left\{
 \varphi:E\to\mathbb R:
 \varphi\ \text{is Borel measurable and}\ 
 \sup_{x\in E}
 \frac{|\varphi(x)|}{1+V(x)}
 <\infty
 \right\}.
\]
From this point onward we assume that $(T_t)$ is defined on $B_V(E)$, or has
its canonical monotone extension to this space, and that $T_tV$ is finite.
This is automatic in the diffusion applications below: if
$\abs{\varphi}\leq C(1+V)$, monotonicity and the Lyapunov estimate imply
$\abs{T_t\varphi}\leq C(1+T_tV)<\infty$.
For $\beta>0$ and $\varphi\in B_V(E)$, define
\begin{align}
 \norm{\varphi}_{\beta,V}
 &:=\sup_{x\in E}\frac{\abs{\varphi(x)}}{1+\beta V(x)},
 \label{eq:weightednorm}\\
 [\varphi]_{\beta,V}
 &:=\sup_{x\neq y}
 \frac{\abs{\varphi(x)-\varphi(y)}}
 {2+\beta V(x)+\beta V(y)}.
 \label{eq:weightedsemi}
\end{align}

The following lemma appeared in \cite{HairerMattingly2011} without a proof. We include a proof for convenience without claiming the originality. 
\begin{lemma}\label{lem:quotient}
For every $\varphi\in B_V(E)$,
\begin{equation}\label{eq:quotient}
 [\varphi]_{\beta,V}
 =\inf_{c\in\R}\norm{\varphi-c}_{\beta,V}.
\end{equation}
In particular, if $[\varphi]_{\beta,V}\leq1$, then for every
$\varepsilon>0$ there is a constant $c$ such that
\begin{equation}\label{eq:normalise}
 \abs{\varphi(x)-c}
 \leq(1+\varepsilon)(1+\beta V(x)),\qquad x\in E.
\end{equation}
\end{lemma}

\begin{proof}
Write $w(x)=1+\beta V(x)$.  If
$\norm{\varphi-c}_{\beta,V}\leq r$, then
\[
 \abs{\varphi(x)-\varphi(y)}
 \leq r(w(x)+w(y)),
\]
and hence $[\varphi]_{\beta,V}\leq r$.

Conversely, suppose $[\varphi]_{\beta,V}<r$.  For every $x$ consider the
closed interval
\[
 I_x=[\varphi(x)-rw(x),\,\varphi(x)+rw(x)].
\]
The definition of the seminorm implies $I_x\cap I_y\neq\varnothing$ for
every pair $x,y$.  By the one-dimensional Helly theorem
\cite{Rockafellar1970}, a family of bounded and closed intervals in $\R$
with pairwise nonempty intersections has a nonempty total intersection. Choose
$c\in\bigcap_x I_x$.  Then
$\abs{\varphi(x)-c}\leq rw(x)$ for every $x$.  Letting
$r\downarrow[\varphi]_{\beta,V}$ proves \eqref{eq:quotient}.
\end{proof}

Let $\mathfrak E_V$ denote the sublinear expectations on $B_V(E)$ satisfying
$\cE[V]<\infty$.  For $\cE_1,\cE_2\in\mathfrak E_V$ define
\begin{equation}\label{eq:metricE}
 D_{\beta,V}(\cE_1,\cE_2)
 :=\sup_{[\varphi]_{\beta,V}\leq1}
 \abs{\cE_1[\varphi]-\cE_2[\varphi]}.
\end{equation}

\begin{lemma}\label{lem:Dfinite}
The quantity in \eqref{eq:metricE} is finite and separates distinct
sublinear expectations in $\mathfrak E_V$.  More precisely,
\begin{equation}\label{eq:Dmoment}
 D_{\beta,V}(\cE_1,\cE_2)
 \leq 2+\beta\bigl(\cE_1[V]+\cE_2[V]\bigr).
\end{equation}
\end{lemma}

\begin{proof}
Fix \(\varphi\in B_V(E)\) with
\[
[\varphi]_{\beta,V}\leq 1,
\]
and let \(\varepsilon>0\). By Lemma~\ref{lem:quotient}, there exists
\(c_\varepsilon\in\mathbb R\) such that, with
\[
\psi:=\varphi-c_\varepsilon,
\]
we have
\begin{equation}\label{eq:Dfinite-normalisation}
 |\psi(x)|
 \leq
 (1+\varepsilon)(1+\beta V(x)),
 \qquad x\in E.
\end{equation}
Translation invariance gives
\[
\cE_i[\varphi]
=
\cE_i[\psi]+c_\varepsilon,
\qquad i=1,2,
\]
and therefore
\[
\bigl|\cE_1[\varphi]-\cE_2[\varphi]\bigr|
=
\bigl|\cE_1[\psi]-\cE_2[\psi]\bigr|.
\]

By \eqref{eq:absbound} and monotonicity,
\[
|\cE_i[\psi]|
\leq
\cE_i[|\psi|]
\leq
(1+\varepsilon)\cE_i[1+\beta V].
\]
Using translation invariance and positive homogeneity,
\[
\cE_i[1+\beta V]
=
1+\beta\cE_i[V].
\]
Consequently,
\begin{align*}
\bigl|\cE_1[\varphi]-\cE_2[\varphi]\bigr|
&\leq
|\cE_1[\psi]|+|\cE_2[\psi]|\\
&\leq
(1+\varepsilon)
\left(
2+\beta\bigl(\cE_1[V]+\cE_2[V]\bigr)
\right).
\end{align*}
Letting \(\varepsilon\downarrow0\) and then taking the supremum over
all \(\varphi\) with \([\varphi]_{\beta,V}\leq1\) yields
\[
D_{\beta,V}(\cE_1,\cE_2)
\leq
2+\beta\bigl(\cE_1[V]+\cE_2[V]\bigr)<\infty.
\]

It remains to prove separation. Suppose that
\[
D_{\beta,V}(\cE_1,\cE_2)=0.
\]
Every \(\varphi\in B_V(E)\) has finite weighted seminorm. Indeed,
\[
[\varphi]_{\beta,V}
\leq
\|\varphi\|_{\beta,V}
\leq
\max\{1,\beta^{-1}\}\|\varphi\|_V
<\infty.
\]
Set \(r=[\varphi]_{\beta,V}\). If \(r=0\), then \(\varphi\) is
constant, so constant preservation gives
\[
\cE_1[\varphi]=\cE_2[\varphi].
\]
If \(r>0\), then
\[
\left[\frac{\varphi}{r}\right]_{\beta,V}=1.
\]
The definition of \(D_{\beta,V}\) therefore gives
\[
\cE_1\left[\frac{\varphi}{r}\right]
=
\cE_2\left[\frac{\varphi}{r}\right].
\]
Positive homogeneity now implies
\[
\cE_1[\varphi]
=
r\,\cE_1\left[\frac{\varphi}{r}\right]
=
r\,\cE_2\left[\frac{\varphi}{r}\right]
=
\cE_2[\varphi].
\]
Thus \(\cE_1=\cE_2\) on \(B_V(E)\), proving that
\(D_{\beta,V}\) separates distinct sublinear expectations.
\end{proof}

The following lemma was proved in \cite{ZhaoZhao2026}. 
\begin{lemma}
\label{lem:E-complete}
The space $(\mathfrak E_V,D_{\beta,V})$ is a complete metric space.
\end{lemma}

\section{From local Harnack to local oscillation contraction}
\label{sec:harnack}

Fix $a>0$ and a nonempty set $C\subset E$. Let $\mathfrak Q$ be a
positive-time space--time set containing the later-time slice
$\{a\}\times C$ and an interior reference point $(s_*,x_*)$ with
$0<s_*<a$.
We start from the following two-sided local Harnack property assumption. 

\renewcommand{\theassumption}{H}
\begin{assumption}
\label{ass:H}
There is $C_H\geq1$ such that the following holds.  For every bounded measurable initial
function $\varphi:E\to \R$, put $u(t,x)=T_t\varphi(x)$ and
\[
 m=\inf_{(t,x)\in\mathfrak Q}u(t,x),\qquad
 M=\sup_{(t,x)\in\mathfrak Q}u(t,x).
\]
The nonnegative functions $v=u-m$ and $w=M-u$ satisfy
\begin{align}
 v(s_*,x_*)&\leq C_H\inf_{x\in C}v(a,x),
 \label{eq:Hupper}\\
 w(s_*,x_*)&\leq C_H\inf_{x\in C}w(a,x).
 \label{eq:Hdual}
\end{align}
\end{assumption}

In applications, \eqref{eq:Hupper} is the parabolic Harnack inequality for
the upper equation, whereas \eqref{eq:Hdual} is Harnack for the conjugate
lower equation.  The next elementary lemma, which reveals the local oscillation contraction from Harnack inequalities, is the local heart of the method.

\begin{lemma}
\label{lem:Hosc}
Under Assumption \ref{ass:H}, 
we have
\begin{equation}\label{eq:localosc}
 \osc_Cu(a,\cdot)
 \leq
 \vartheta\,\osc_{\mathfrak Q}u,
 \qquad
 \vartheta:=1-\frac{1}{C_H}<1.
\end{equation}
In particular, if $u(t,x)=T_t\varphi(x)$, then
\[
 \osc_C(T_a\varphi)
 \leq
 \vartheta\,
 \osc_{\mathfrak Q}(T_\cdot\varphi).
\]
\end{lemma}

\begin{proof}
Put
\[
 A:=M-m=\osc_{\mathfrak Q}u.
\]
Since $v+w=M-m$ throughout $\mathfrak Q$, at the reference point
\[
 v(s_*,x_*)+w(s_*,x_*)=A.
\]
The two assumed Harnack estimates imply
\[
 \inf_{x\in C}v(a,x)
 \geq
 \frac{v(s_*,x_*)}{C_H},
 \qquad
 \inf_{x\in C}w(a,x)
 \geq
 \frac{w(s_*,x_*)}{C_H}.
\]
Adding these inequalities gives
\begin{equation}\label{eq:Hosc-sum}
 \inf_{x\in C}v(a,x)
 +
 \inf_{x\in C}w(a,x)
 \geq
 \frac{A}{C_H}.
\end{equation}
Moreover,
\[
 \inf_{x\in C}u(a,x)
 =
 m+\inf_{x\in C}v(a,x)
\]
and
\[
 \sup_{x\in C}u(a,x)
 =
 M-\inf_{x\in C}w(a,x).
\]
Consequently, by \eqref{eq:Hosc-sum},
\begin{align*}
 \osc_Cu(a,\cdot)
 &=
 A
 -\inf_{x\in C}v(a,x)
 -\inf_{x\in C}w(a,x)\\
 &\leq
 \left(1-\frac{1}{C_H}\right)A\\
 &=
 \vartheta\,\osc_{\mathfrak Q}u.
\end{align*}
This proves \eqref{eq:localosc}.
\end{proof}

\begin{remark}\label{rem:both}
Harnack for $u-m$ alone only improves the lower side of the range when
$u(s_*,x_*)$ is sufficiently far above $m$.  If the reference value lies
near $m$, the useful quantity is $M-u$, which evolves according to the lower,
not the upper, nonlinear equation. Therefore, both estimates in the proof of
Lemma~\ref{lem:Hosc} are essential.
\end{remark}

\section{The abstract weighted contraction theorem}
\label{sec:abstract}

Suppose that $V:E\to[0,\infty)$ has bounded sublevel sets
\[
 C_R:=\{x\in E:V(x)\leq R\}.
\]
We impose two assumptions at a fixed time $a>0$.

\renewcommand{\theassumption}{LD}
\begin{assumption}[Lyapunov drift]
\label{ass:L}
There exist $a>0$, $\gamma\in(0,1)$ and $K\geq0$ such that
\begin{equation}\label{eq:Lyap-a}
 T_aV\leq\gamma V+K.
\end{equation}
Moreover, there exist constants $\Gamma_a<\infty$ and
$K_a<\infty$ such that
\begin{equation}\label{eq:Lyap-short}
 T_rV\leq\Gamma_aV+K_a,
 \qquad 0\leq r\leq a.
\end{equation}
\end{assumption}

Assumption~\ref{ass:L} also gives control over every finite time
interval. Indeed, if $t=na+r$, where $n\in\mathbb N_0$ and
$r\in[0,a)$, iteration of \eqref{eq:Lyap-a} gives
\[
 T_{na}V
 \leq
 \gamma^nV
 +
 K\sum_{j=0}^{n-1}\gamma^j
 \leq
 \gamma^nV+\frac{K}{1-\gamma}.
\]
Using the semigroup property, monotonicity, positive homogeneity and
constant preservation, we obtain
\begin{align}
 T_tV
 &=
 T_r(T_{na}V) \notag\\
 &\leq
 \gamma^nT_rV+\frac{K}{1-\gamma} \notag\\
 &\leq
 \Gamma_aV+K_a+\frac{K}{1-\gamma}.
 \label{eq:Lyap-finite-horizon}
\end{align}
In particular, for every $T<\infty$, there exist
$\Gamma_T,K_T<\infty$ such that
\begin{equation}\label{eq:Lyap-general-horizon}
 T_tV\leq\Gamma_TV+K_T,
 \qquad 0\leq t\leq T.
\end{equation}

Choose
\begin{equation}\label{eq:Rchoice}
 R>\frac{2K}{1-\gamma}.
\end{equation}

\renewcommand{\theassumption}{LO}
\begin{assumption}[Local oscillation contraction]
\label{ass:LOC}
There exist constants
\[
 T_R\in[a,\infty),
 \qquad
 \widehat R\in[R,\infty),
 \qquad
 \vartheta_R\in[0,1),
\]
and a positive-time space--time set $\mathfrak Q_R$ satisfying
\begin{equation}\label{eq:QRgeometry}
 \{a\}\times C_R
 \subset
 \mathfrak Q_R
 \subset
 (0,T_R]\times C_{\widehat R},
\end{equation}
such that, for every $\varphi\in B_V(E)$, with
\[
 u(t,x):=T_t\varphi(x),
\]
one has
\begin{equation}\label{eq:abstract-local-osc}
 \osc_{C_R}(T_a\varphi)
 \leq
 \vartheta_R\osc_{\mathfrak Q_R}u.
\end{equation}
Here
\[
 \osc_{\mathfrak Q_R}u
 :=
 \sup_{(t,x)\in\mathfrak Q_R}u(t,x)
 -
 \inf_{(t,x)\in\mathfrak Q_R}u(t,x).
\]
\end{assumption}

The local Lyapunov bound needed below is now automatic. Indeed,
\eqref{eq:Lyap-general-horizon} and \eqref{eq:QRgeometry} give
\begin{equation}\label{eq:LR}
 L_R
 :=
 \sup_{(t,x)\in\mathfrak Q_R}T_tV(x)
 \leq
 \Gamma_{T_R}\widehat R+K_{T_R}
 <\infty.
\end{equation}

\begin{remark}
\label{rem:LOC-Harnack}
Assumption~\ref{ass:LOC} is stated directly in terms of the
oscillation inequality used in the proof below. Neither a reference
point nor a Harnack inequality is part of the abstract assumption.

When Assumption~\ref{ass:LOC} is verified through the two-sided
Harnack property of Assumption~\ref{ass:H} and
Lemma~\ref{lem:Hosc}, the set $\mathfrak Q_R$ is chosen to contain
the corresponding finite Harnack chain and an interior reference
point $(s_*,x_*)$. A finite Harnack chain has bounded time
projection, and in the diffusion applications its spatial projection
is contained in a bounded $V$-sublevel set. Thus constants
$T_R,\widehat R<\infty$ satisfying \eqref{eq:QRgeometry} can be
chosen. The reference point is used only to deduce
\eqref{eq:abstract-local-osc}; it plays no role in the weighted
contraction argument.
\end{remark}

The following weighted nonlinear Harris--oscillation contraction
theorem is the principal abstract result of this paper. It replaces
the usual local minorisation condition by a direct contraction of
oscillations and is therefore well suited to nonlinear Markov
semigroups.

\begin{theorem}
\label{thm:weighted}
Under Assumptions~\ref{ass:L} and \ref{ass:LOC}, there exist
$\beta>0$ and $q\in(0,1)$ such that
\begin{equation}\label{eq:maincontraction}
 [T_a\varphi]_{\beta,V}
 \leq
 q[\varphi]_{\beta,V},
 \qquad
 \varphi\in B_V(E).
\end{equation}
More precisely, choose $\beta>0$ sufficiently small that
\begin{equation}\label{eq:qin}
 q_{\mathrm{in}}
 :=
 \vartheta_R(1+\beta L_R)
 <1,
\end{equation}
where $L_R$ is given by \eqref{eq:LR}, and set
\begin{align}
 \gamma_0
 &:=
 \gamma+\frac{2K}{R}<1,
 \label{eq:gamma0}\\
 q_{\mathrm{out}}
 &:=
 \frac{2+\beta\gamma_0R}
      {2+\beta R}
 <1,
 \label{eq:qout}\\
 q
 &:=
 q_{\mathrm{in}}\vee q_{\mathrm{out}}.
\end{align}
\end{theorem}

\begin{proof}
By positive homogeneity, it is enough to prove the result when
\[
 [\varphi]_{\beta,V}\leq1.
\]
By Lemma~\ref{lem:quotient}, translation invariance and a limiting
argument, we may assume that
\begin{equation}\label{eq:phinormproof}
 |\varphi(x)|\leq1+\beta V(x),
 \qquad x\in E.
\end{equation}
The order properties of $T_t$ give
\begin{equation}\label{eq:TphiV}
 |T_t\varphi(x)|
 \leq
 1+\beta T_tV(x).
\end{equation}

Fix $x,y\in E$ and write
\[
 S:=V(x)+V(y).
\]
We consider separately the cases $S\leq R$ and $S\geq R$.

Suppose first that $S\leq R$. Then $x,y\in C_R$. By
\eqref{eq:TphiV} and \eqref{eq:LR},
\[
 \osc_{\mathfrak Q_R}(T_\cdot\varphi)
 \leq
 2(1+\beta L_R).
\]
Assumption~\ref{ass:LOC} therefore yields
\[
 |T_a\varphi(x)-T_a\varphi(y)|
 \leq
 2\vartheta_R(1+\beta L_R).
\]
Consequently,
\begin{align*}
 \frac{|T_a\varphi(x)-T_a\varphi(y)|}
      {2+\beta V(x)+\beta V(y)}
 &\leq
 \frac{2\vartheta_R(1+\beta L_R)}
      {2+\beta S}\\
 &\leq
 \vartheta_R(1+\beta L_R)
 =
 q_{\mathrm{in}}.
\end{align*}

Suppose next that $S\geq R$. Using \eqref{eq:TphiV} and
\eqref{eq:Lyap-a}, we obtain
\begin{align*}
 |T_a\varphi(x)-T_a\varphi(y)|
 &\leq
 2+\beta\bigl(T_aV(x)+T_aV(y)\bigr)\\
 &\leq
 2+\beta\gamma S+2\beta K\\
 &\leq
 2+\beta\gamma_0S,
\end{align*}
where the last inequality follows from $S\geq R$ and
\[
 \gamma_0=\gamma+\frac{2K}{R}.
\]
The choice \eqref{eq:Rchoice} ensures that $\gamma_0<1$. Hence the
function
\[
 s\longmapsto
 \frac{2+\beta\gamma_0s}{2+\beta s}
\]
is decreasing on $[0,\infty)$. It follows that
\begin{align*}
 \frac{|T_a\varphi(x)-T_a\varphi(y)|}
      {2+\beta V(x)+\beta V(y)}
 &\leq
 \frac{2+\beta\gamma_0S}{2+\beta S}\\
 &\leq
 \frac{2+\beta\gamma_0R}{2+\beta R}
 =
 q_{\mathrm{out}}.
\end{align*}

Taking the supremum over $x\neq y$ and combining the two cases gives
\[
 [T_a\varphi]_{\beta,V}
 \leq
 \bigl(q_{\mathrm{in}}\vee q_{\mathrm{out}}\bigr)
 [\varphi]_{\beta,V}.
\]
Since both $q_{\mathrm{in}}$ and $q_{\mathrm{out}}$ are strictly
smaller than one, this proves \eqref{eq:maincontraction}.
\end{proof}

\begin{remark}
The abstract theorem may be viewed as an operator-level nonlinear Harris theorem. In the classical Harris framework, local mixing is commonly expressed through a Doeblin minorisation or small-set condition. Here, it is replaced by a direct local oscillation-contraction assumption. Every suitable Doeblin minorisation implies such an oscillation contraction, but the latter does not require the existence of a common minorising measure. This weaker formulation is particularly useful for nonlinear Markov semigroups and controlled G-diffusions, for which local contraction may be established analytically through upper and dual Harnack inequalities even when a common kernel minorisation is difficult to formulate or verify.
\end{remark}

\begin{remark}
Theorem~\ref{thm:weighted} is formulated entirely at the semigroup
level. Its local assumption is the oscillation contraction estimate
\[
 \osc_{C_R}(T_a\varphi)
 \leq
 \vartheta_R\,
 \osc_{\mathfrak Q_R}(T_\cdot\varphi),
 \qquad \vartheta_R<1,
\]
rather than a particular probabilistic minorisation or analytic
Harnack condition. Consequently, the theorem is not restricted to
diffusion semigroups or to a specific representation of the nonlinear
expectation.

A two-sided parabolic Harnack inequality is one convenient way of
verifying the local oscillation contraction, but it is not part of the
abstract theorem itself. In other applications the same estimate may
follow from a coupling argument \cite{ChenLi1989, Eberle2016, HairerMattinglyScheutzow2011}, a common minorisation, a strong
Feller--irreducibility mechanism, controllability, or another local
regularisation estimate. See Zhao--Zhao \cite{ZhaoZhao2026} for a nonlinear coupling argument. The abstract result therefore separates the
global recurrence mechanism, expressed by the Lyapunov condition,
from the local mixing mechanism, expressed by oscillation contraction.
\end{remark}

\begin{remark}
\label{rem:Doeblin-oscillation}
The local oscillation-contraction condition is weaker than the usual
Doeblin minorisation.  Indeed, let $P_a$ be a linear Markov operator and
suppose that there exist $\varepsilon\in(0,1]$ and a probability measure
$\nu$ such that
\begin{equation}\label{eq:Doeblin}
P_a(x,\cdot)\geq\varepsilon\nu(\cdot),
\qquad x\in C_R.
\end{equation}
For every $x\in C_R$, there is then a probability measure
$\widehat P_a(x,\cdot)$ such that
\[
P_a(x,\cdot)
=
\varepsilon\nu(\cdot)
+
(1-\varepsilon)\widehat P_a(x,\cdot).
\]
Consequently, for every $\varphi\in B_b(E)$,
\[
P_a\varphi(x)
=
\varepsilon\nu(\varphi)
+
(1-\varepsilon)\widehat P_a\varphi(x).
\]
The first term is independent of $x$.  Hence, for $x,y\in C_R$,
\begin{align*}
\abs{P_a\varphi(x)-P_a\varphi(y)}
&=
(1-\varepsilon)
\abs{
\widehat P_a\varphi(x)-\widehat P_a\varphi(y)
}\\
&\leq
(1-\varepsilon)\osc_E(\varphi).
\end{align*}
Taking the supremum over $x,y\in C_R$ gives
\begin{equation}\label{eq:Doeblin-oscillation}
\osc_{C_R}(P_a\varphi)
\leq
(1-\varepsilon)\osc_E(\varphi).
\end{equation}

More generally, suppose that a sublinear Markov operator has the
representation
\[
T_a\varphi(x)
=
\sup_{\alpha\in\mathcal A}
P_a^\alpha\varphi(x),
\]
and that the common minorisation
\[
P_a^\alpha(x,\cdot)
\geq
\varepsilon\nu(\cdot)
\]
holds for every $x\in C_R$ and every $\alpha\in\mathcal A$.  Writing
\[
P_a^\alpha(x,\cdot)
=
\varepsilon\nu(\cdot)
+
(1-\varepsilon)\widehat P_a^\alpha(x,\cdot)
\]
gives
\[
T_a\varphi(x)
=
\varepsilon\nu(\varphi)
+
(1-\varepsilon)
\sup_{\alpha\in\mathcal A}
\widehat P_a^\alpha\varphi(x).
\]
Therefore,
\[
\osc_{C_R}(T_a\varphi)
\leq
(1-\varepsilon)\osc_E(\varphi).
\]
Thus a common Doeblin minorisation implies an operator-level
oscillation contraction. The converse fails even for a linear Markov
operator on a finite state space. Indeed, let
\[
 E=C_R=\{1,2,3\}
\]
and define
\[
 \begin{aligned}
 P_a(1,\cdot)&=\tfrac12(\delta_1+\delta_2),\\
 P_a(2,\cdot)&=\tfrac12(\delta_2+\delta_3),\\
 P_a(3,\cdot)&=\tfrac12(\delta_1+\delta_3).
 \end{aligned}
\]
For every function $\varphi:E\to\R$, the three pairwise differences of
$P_a\varphi$ are
\[
 \frac12\bigl(\varphi(1)-\varphi(3)\bigr),\qquad
 \frac12\bigl(\varphi(2)-\varphi(3)\bigr),\qquad
 \frac12\bigl(\varphi(2)-\varphi(1)\bigr),
\]
up to their signs. Hence
\[
 \osc_E(P_a\varphi)
 =\frac12\osc_E(\varphi).
\]
Nevertheless, no $\varepsilon>0$ and probability measure $\nu$ can satisfy
\[
 P_a(x,\cdot)\geq\varepsilon\nu(\cdot),
 \qquad x\in E.
\]
Indeed, the zero masses of $P_a(1,\cdot)$ at $3$, of
$P_a(2,\cdot)$ at $1$, and of $P_a(3,\cdot)$ at $2$ would imply
\[
 \nu(\{3\})=\nu(\{1\})=\nu(\{2\})=0,
\]
which is impossible for a probability measure.

Thus, already in the linear setting, oscillation contraction expresses a
uniform pairwise overlap of transition probabilities, whereas a Doeblin
minorisation requires a single component common to all of them. The latter
is strictly stronger. Theorem~\ref{thm:weighted} takes the weaker
oscillation contraction itself as the structural hypothesis. This is
particularly useful for nonlinear Markov semigroups, for which a single
reference measure $\nu$ may be unavailable or unnatural, whereas analytic
estimates for both the upper equation and its dual can still yield the
required local contraction.
\end{remark}

\begin{remark}
When \(T_a=P\) is a linear Markov operator, the seminorm
\[
 [\varphi]_{\beta,V}
 =
 \inf_{c\in\R}\|\varphi-c\|_{\beta,V}
\]
is the quotient form of the weighted supremum norm, and its dual on
probability measures is the corresponding weighted total-variation
distance. In this case Theorem~\ref{thm:weighted} is a variant of the
weighted Harris contraction theorem of Hairer and Mattingly
\cite{HairerMattingly2011}; see also Meyn and Tweedie
\cite{MeynTweedie2009}.

The distinction is that Assumption~\ref{ass:LOC} postulates directly a
local contraction of oscillations rather than a common Doeblin
minorisation. In the linear setting this is a local Dobrushin-type
condition controlling the pairwise overlap of transition probabilities.
As the three-state example in Remark~\ref{rem:Doeblin-oscillation} shows,
it is strictly weaker than requiring one measure to minorise all of the
transition probabilities simultaneously. The principal purpose of the
present formulation is that it continues to make sense for sublinear Markov
semigroups, where there may be no single transition kernel or common
minorising measure.

\end{remark}

\begin{remark}
\label{rem:long-run-risk}
Let the test function \(\varphi\) be interpreted as a loss.  The
invariant expectation \(\mathcal E_\infty\) defines the coherent risk
functional
\[
 \rho_\infty(\varphi):=\mathcal E_\infty(\varphi).
\]
Its invariance can be written as
\[
 \rho_\infty(T_t\varphi)=\rho_\infty(\varphi),
 \qquad t\geq0.
\]
Thus \(\rho_\infty\) evaluates a loss in a manner that is unchanged by
the \(G\)-diffusion dynamics.  In this sense it is the stationary
long-run risk evaluation associated with the uncertain diffusion.

Moreover, the exponential convergence established above has the
following risk interpretation: for every admissible loss
\(\varphi\), the finite-horizon worst-case risk
\[
 x\longmapsto T_t\varphi(x)
\]
forgets the initial state exponentially fast and converges to the
state-independent quantity \(\rho_\infty(\varphi)\).  Consequently,
the invariant expectation describes not only an invariant object for
the nonlinear semigroup, but also the asymptotic risk profile of the
system.

Although the present paper treats autonomous semigroups, the same
interpretation points toward a nonautonomous theory. In that setting,
the asymptotic risk object is not a single stationary functional but
an entrance family
\[
 (\rho_t)_{t\in\mathbb R}.
\]
More precisely, if $(T_{s,t})_{s\leq t}$ is a nonlinear evolution and
$(\mathcal E_t)_{t\in\mathbb R}$ is an entrance expectation, define
\[
 \rho_t(\varphi):=\mathcal E_t(\varphi).
\]
The entrance identity
\[
 \mathcal E_s(T_{s,t}\varphi)=\mathcal E_t(\varphi),
 \qquad s\leq t,
\]
then becomes
\[
 \rho_s(T_{s,t}\varphi)=\rho_t(\varphi),
 \qquad s\leq t.
\]
This is the natural entrance-consistency relation for a nonstationary
family of coherent long-run risk evaluations and is the nonautonomous
counterpart of
\[
 \rho_\infty(T_t\varphi)=\rho_\infty(\varphi).
\]

If \(\varphi\) represents a profit rather than a loss, one should use
\(\rho_t^{\,\mathrm{profit}}(\varphi)
 :=\mathcal E_t(-\varphi)\), in accordance with the usual sign
convention for monetary risk measures.
\end{remark}

%
The following proposition reveals 
operator-level convergence in the space of sublinear expectation operators.
\begin{proposition}
\label{prop:operator-convergence}
Let
\[
\mathcal P\mathcal E:=\mathcal E\circ T_a,
\qquad
\mathcal E\in\mathfrak E_V.
\]
Under the assumptions of Theorem~\ref{thm:weighted},
$\mathcal P$ maps $\mathfrak E_V$ into itself and is a strict
contraction:
\begin{equation}
\label{zhao10a}
D_{\beta,V}(\mathcal P\mathcal E_1,\mathcal P\mathcal E_2)
\leq
qD_{\beta,V}(\mathcal E_1,\mathcal E_2),
\qquad q\in(0,1).
\end{equation}
Consequently, there exists a unique
$\widetilde{\mathbb E}\in\mathfrak E_V$ satisfying
\[
\mathcal P\widetilde{\mathbb E}
=
\widetilde{\mathbb E},
\]
or equivalently,
\[
\widetilde{\mathbb E}[T_a\varphi]
=
\widetilde{\mathbb E}[\varphi],
\qquad
\varphi\in B_V(E).
\]
Moreover, for every $\mathcal E_0\in\mathfrak E_V$ and every
$n\geq0$,
\begin{equation}\label{eq:operator-rate}
D_{\beta,V}
\bigl(
\mathcal E_0\circ T_{na},
\widetilde{\mathbb E}
\bigr)
\leq
q^nD_{\beta,V}
\bigl(
\mathcal E_0,
\widetilde{\mathbb E}
\bigr).
\end{equation}
One also has the a posteriori estimate
\begin{equation}\label{eq:operator-rate-aposteriori}
D_{\beta,V}
\bigl(
\mathcal E_0\circ T_{na},
\widetilde{\mathbb E}
\bigr)
\leq
\frac{q^n}{1-q}
D_{\beta,V}
\bigl(
\mathcal E_0\circ T_a,
\mathcal E_0
\bigr).
\end{equation}
\end{proposition}

\begin{proof}
The Lyapunov condition gives
\[
(\mathcal P\mathcal E)[V]
=
\mathcal E[T_aV]
\leq
\gamma\mathcal E[V]+K<\infty.
\]
Thus $\mathcal P\mathcal E\in\mathfrak E_V$ whenever
$\mathcal E\in\mathfrak E_V$.

We first prove (\ref{zhao10a}). Indeed, if $[\varphi]_{\beta,V}\le 1$, then Theorem~\ref{thm:weighted}
yields
\[
[T_a\varphi]_{\beta,V}\le q.
\]
Hence, by the positive homogeneity of $[\cdot]_{\beta,V}$ and the definition of
$D_{\beta,V}$,
\[
\begin{aligned}
\bigl|\mathcal E_1[T_a\varphi]-\mathcal E_2[T_a\varphi]\bigr|
&=
q\left|
\mathcal E_1\!\left[q^{-1}T_a\varphi\right]
-\mathcal E_2\!\left[q^{-1}T_a\varphi\right]
\right|  \\
&\le q\,D_{\beta,V}(\mathcal E_1,\mathcal E_2),
\end{aligned}
\]
because
\[
\left[q^{-1}T_a\varphi\right]_{\beta,V}\le 1.
\]
Taking the supremum over all $\varphi$ satisfying
$[\varphi]_{\beta,V}\le 1$ therefore gives the desired result (\ref{zhao10a}).

Thus, $\mathcal P$ is a strict
contraction on the complete metric space
$(\mathfrak E_V,D_{\beta,V})$.  Banach's fixed-point theorem therefore
gives a unique fixed point $\widetilde{\mathbb E}\in\mathfrak E_V$.

Since
\[
\mathcal P^n\mathcal E_0
=
\mathcal E_0\circ T_{na},
\]
iteration of the contraction inequality gives
\[
D_{\beta,V}
\bigl(
\mathcal E_0\circ T_{na},
\widetilde{\mathbb E}
\bigr)
=
D_{\beta,V}
\bigl(
\mathcal P^n\mathcal E_0,
\mathcal P^n\widetilde{\mathbb E}
\bigr)
\leq
q^nD_{\beta,V}
\bigl(
\mathcal E_0,
\widetilde{\mathbb E}
\bigr).
\]
This proves \eqref{eq:operator-rate}.

For the second estimate, observe that
\begin{align*}
D_{\beta,V}
\bigl(
\mathcal P^n\mathcal E_0,
\widetilde{\mathbb E}
\bigr)
&\leq
\sum_{k=n}^{\infty}
D_{\beta,V}
\bigl(
\mathcal P^{k+1}\mathcal E_0,
\mathcal P^k\mathcal E_0
\bigr)\\
&\leq
\sum_{k=n}^{\infty}
q^k
D_{\beta,V}
\bigl(
\mathcal P\mathcal E_0,
\mathcal E_0
\bigr)\\
&=
\frac{q^n}{1-q}
D_{\beta,V}
\bigl(
\mathcal E_0\circ T_a,
\mathcal E_0
\bigr).
\end{align*}
\end{proof}

The following result is the convergence in the sense of point evaluations. 
\begin{corollary}
\label{cor:point-evaluation-convergence}
Under the assumptions of Theorem~\ref{thm:weighted}, let
\(\widetilde{\mathbb E}\) be the unique fixed point of
$
\mathcal P: \mathfrak E_V \to \mathfrak E_V$. Then, for every \(x\in E\) and \(n\ge0\),
\begin{equation}
\label{eq:point-operator-rate}
D_{\beta,V}
\bigl(
\delta_x\circ T_{na},
\widetilde{\mathbb E}
\bigr)
\le
\frac{q^n}{1-q}
\left[
2+\beta\bigl((1+\gamma)V(x)+K\bigr)
\right],
\end{equation}
where
\[
\delta_x[\varphi]:=\varphi(x).
\]
Consequently, for every \(\varphi\in B_V(E)\),
\begin{equation}
\label{eq:pointwise-discrete-rate}
\left|
T_{na}\varphi(x)-\widetilde{\mathbb E}[\varphi]
\right|
\le
\frac{q^n}{1-q}
\left[
2+\beta\bigl((1+\gamma)V(x)+K\bigr)
\right]
[\varphi]_{\beta,V}.
\end{equation}
\end{corollary}

\begin{proof}
Since \(V(x)<\infty\), the point-evaluation functional
\[
\delta_x[\varphi]=\varphi(x)
\]
belongs to \(\mathfrak E_V\). Applying the a posteriori estimate from
Proposition~\ref{prop:operator-convergence} with
\(\mathcal E_0=\delta_x\), we obtain
\begin{equation}
\label{eq:point-posteriori-bound}
D_{\beta,V}
\bigl(
\delta_x\circ T_{na},
\widetilde{\mathbb E}
\bigr)
\le
\frac{q^n}{1-q}
D_{\beta,V}
\bigl(
\delta_x\circ T_a,
\delta_x
\bigr).
\end{equation}

It therefore remains to estimate the one-step distance on the
right-hand side. Let \(\varphi\in B_V(E)\) satisfy
\[
[\varphi]_{\beta,V}\le1.
\]
By Lemma~\ref{lem:quotient}, for every
\(\varepsilon>0\) there exists \(c_\varepsilon\in\mathbb R\) such that
\[
|\varphi(y)-c_\varepsilon|
\le
(1+\varepsilon)\bigl(1+\beta V(y)\bigr),
\qquad y\in E.
\]
Translation invariance of \(T_a\) gives
\[
T_a\varphi(x)-\varphi(x)
=
T_a(\varphi-c_\varepsilon)(x)
-
\bigl(\varphi(x)-c_\varepsilon\bigr).
\]
Moreover, monotonicity and sublinearity imply
\[
\bigl|T_a(\varphi-c_\varepsilon)(x)\bigr|
\le
T_a|\varphi-c_\varepsilon|(x).
\]
It follows that
\[
\begin{aligned}
|T_a\varphi(x)-\varphi(x)|
&\le
T_a|\varphi-c_\varepsilon|(x)
+
|\varphi(x)-c_\varepsilon|\\
&\le
(1+\varepsilon)
\left[
2+\beta\bigl(T_aV(x)+V(x)\bigr)
\right].
\end{aligned}
\]
Using the Lyapunov estimate
\[
T_aV(x)\le\gamma V(x)+K,
\]
we conclude that
\[
|T_a\varphi(x)-\varphi(x)|
\le
(1+\varepsilon)
\left[
2+\beta\bigl((1+\gamma)V(x)+K\bigr)
\right].
\]
Letting \(\varepsilon\downarrow0\) and taking the supremum over all
\(\varphi\) with \([\varphi]_{\beta,V}\le1\) gives
\begin{equation}
\label{eq:point-one-step-bound}
D_{\beta,V}
\bigl(
\delta_x\circ T_a,
\delta_x
\bigr)
\le
2+\beta\bigl((1+\gamma)V(x)+K\bigr).
\end{equation}
Combining \eqref{eq:point-posteriori-bound} and
\eqref{eq:point-one-step-bound} proves
\eqref{eq:point-operator-rate}.

Finally, suppose that \([\varphi]_{\beta,V}>0\). Applying the
definition of \(D_{\beta,V}\) to
\(
\frac{\varphi}{[\varphi]_{\beta,V}}
\)
yields
\[
\begin{aligned}
\left|
T_{na}\varphi(x)-\widetilde{\mathbb E}[\varphi]
\right|
&=
\left|
(\delta_x\circ T_{na})[\varphi]
-
\widetilde{\mathbb E}[\varphi]
\right|\\
&\le
[\varphi]_{\beta,V}
D_{\beta,V}
\bigl(
\delta_x\circ T_{na},
\widetilde{\mathbb E}
\bigr).
\end{aligned}
\]
Estimate \eqref{eq:pointwise-discrete-rate} now follows from
\eqref{eq:point-operator-rate}. If
\([\varphi]_{\beta,V}=0\), then \(\varphi\) is constant, and both sides
of \eqref{eq:pointwise-discrete-rate} vanish.
\end{proof}

The following is the exponential convergence to the invariant expectation.
\begin{corollary}
\label{thm:invariant}
Under the assumptions of Theorem~\ref{thm:weighted}, for every bounded $\varphi$ there are $C,\rho>0$, independent of
$\varphi$ and $x$, such that
\begin{equation}\label{eq:continuousrate}
 \abs{T_t\varphi(x)-\Et[\varphi]}
 \leq C e^{-\rho t}(1+V(x))\osc_E(\varphi),
 \qquad t\geq0.
\end{equation}
Consequently, the convergence is uniform on every set on which $V$ is
bounded.
\end{corollary}

\begin{proof} By Proposition~\ref{prop:operator-convergence}, \(
T_a
\)
has a unique fixed point in $\mathfrak E_V$, which we denote by $\widetilde{\mathbb E}$.
We first show that $\widetilde{\mathbb E}$ is invariant under the whole
semigroup. Fix $0\leq r<a$. By the short-time Lyapunov estimate
\eqref{eq:Lyap-short},
\[
\widetilde{\mathbb E}\circ T_r\in\mathfrak E_V.
\]
Moreover, the semigroup property gives
\begin{align*}
(\widetilde{\mathbb E}\circ T_r)\circ T_a
&=\widetilde{\mathbb E}\circ T_{r+a}\\
&=(\widetilde{\mathbb E}\circ T_a)\circ T_r\\
&=\widetilde{\mathbb E}\circ T_r.
\end{align*}
Hence $\widetilde{\mathbb E}\circ T_r$ is another fixed point of
$\mathcal P$. By uniqueness of the fixed point in $\mathfrak E_V$,
\[
\widetilde{\mathbb E}\circ T_r
=
\widetilde{\mathbb E},
\qquad 0\leq r<a.
\]
Combining this with the $T_a$-invariance and the semigroup property
yields
\begin{equation}\label{eq:invalltime}
\widetilde{\mathbb E}\circ T_t
=
\widetilde{\mathbb E},
\qquad t\geq0.
\end{equation}
In particular, any sublinear expectation in $\mathfrak E_V$ that is
invariant under every $T_t$ is $T_a$-invariant and must therefore equal
$\widetilde{\mathbb E}$.

Now write
\[
t=na+r,
\qquad n\geq0,\quad 0\leq r<a.
\]
Applying \eqref{eq:pointwise-discrete-rate} to $T_r\varphi$ and using
\eqref{eq:invalltime}, we obtain
\begin{align*}
\abs{T_t\varphi(x)-\widetilde{\mathbb E}[\varphi]}
&=
\abs{
T_{na}(T_r\varphi)(x)
-\widetilde{\mathbb E}[T_r\varphi]
}\\
&\leq
\frac{q^n}{1-q}
\left[
2+\beta\bigl((1+\gamma)V(x)+K\bigr)
\right]
[T_r\varphi]_{\beta,V}.
\end{align*}
If $\varphi$ is bounded, then the oscillation non-expansion property
implies
\[
[T_r\varphi]_{\beta,V}
\leq
\frac12\osc_E(T_r\varphi)
\leq
\frac12\osc_E(\varphi).
\]
Furthermore, there is a structural constant $C_0>0$ such that
\[
2+\beta\bigl((1+\gamma)V(x)+K\bigr)
\leq C_0(1+V(x)),
\qquad x\in E.
\]
For example, one may take
\[
C_0=\max\{2+\beta K,\,\beta(1+\gamma)\}.
\]

Set
\[
\rho:=-\frac{\log q}{a}>0.
\]
Since \(0\leq r<a\),
\[
q^n
=e^{-\rho na}
=e^{-\rho(t-r)}
\leq e^{\rho a}e^{-\rho t}
=q^{-1}e^{-\rho t}.
\]
Consequently,
\[
\abs{T_t\varphi(x)-\widetilde{\mathbb E}[\varphi]}
\leq
C e^{-\rho t}(1+V(x))\osc_E(\varphi),
\]
where, for example,
\[
C=\frac{C_0}{2q(1-q)}.
\]
The constants $C$ and $\rho$ are independent of $\varphi$ and $x$.

Finally, if $A\subset E$ and
\[
\sup_{x\in A}V(x)<\infty,
\]
taking the supremum over $x\in A$ proves that
$T_t\varphi(x)\to\widetilde{\mathbb E}[\varphi]$ uniformly on $A$.
\end{proof}

\begin{remark}
\label{rem:class}
On an unbounded state space, weighted contraction gives uniqueness in
$\mathfrak E_V$, equivalently among invariant expectations having finite
Lyapunov moment.  If an expectation is given only as an abstract functional
on $B_b(E)$, the statement $\cE[V]<\infty$ is not even defined.  Uniqueness
among all such abstract functionals requires an additional regularity or
tightness hypothesis permitting extension to $V$.  This distinction
disappears on a compact state space, where $V$ may be chosen bounded. See Feng and Zhao \cite{FengZhao2026} for the smooth manifold without boundary case. 
\end{remark}

\section{A linear benchmark: ordinary uniformly elliptic SDEs}
\label{sec:classical}

The invariant-measure and exponential-mixing conclusions in this
section are classical consequences of Harris-type ergodic theory.
For diffusions, the prototype is the Lyapunov-function
criterion of Khas'minskii \cite{Khasminskii1960}, which already yields
the existence of a unique invariant measure and the stabilisation of
the associated parabolic Cauchy problem; see also Bhattacharya
\cite{Bhattacharya1978} for recurrence criteria and the existence of
invariant measures in the multidimensional case. The general Harris-type
theory is given in Meyn and Tweedie \cite{MeynTweedie2009}, and its
weighted-norm formulation in Hairer and Mattingly
\cite{HairerMattingly2011}. We include the argument for
two reasons. First, it shows that the abstract oscillation-contraction
theorem recovers the standard weighted exponential ergodicity result
for uniformly elliptic diffusions. Second, it isolates in the linear
setting the Harnack mechanism that will later be applied separately to
the upper and conjugate lower Bellman equations. No novelty is claimed
for the invariant-measure conclusion itself.

Consider the It\^o SDE on $\R^d$
\begin{equation}\label{eq:classicalSDE}
 dX_t=b(X_t)\,dt+\sigma(X_t)\,dW_t,
 \qquad X_0=x,
\end{equation}
where $W$ is an $m$-dimensional Brownian motion and
$\sigma:\R^d\to\R^{d\times m}$. Let
\[
 P_t\varphi(x)=E[\varphi(X_t^x)].
\]

\renewcommand{\theassumption}{C}
\begin{assumption}[Classical diffusion hypotheses]\label{ass:SDE}
\ 
\begin{enumerate}[label=\textup{(S\arabic*)}]
 \item $b$ and $\sigma$ are globally Lipschitz and have at most linear
 growth. The functions themselves need not be bounded.

 \item With $V(x)=\abs{x}^2$, there are $c>0$ and $C_0\geq0$ such that
 \begin{equation}\label{eq:classicalLyap}
  2\langle x,b(x)\rangle
  +\norm{\sigma(x)}_{\mathrm{HS}}^2
  \leq
  -c\abs{x}^2+C_0.
 \end{equation}

 \item On every bounded set $O$ needed below there are
 $0<\lambda_O\leq\Lambda_O<\infty$ such that
 \begin{equation}\label{eq:classicalell}
  \lambda_O I_d
  \leq
  \sigma(x)\sigma(x)^{\mathsf T}
  \leq
  \Lambda_O I_d,
  \qquad x\in O.
 \end{equation}
\end{enumerate}
\end{assumption}

The upper bound in \eqref{eq:classicalell} is automatic on bounded
sets under continuity of $\sigma$. The lower bound follows if
$\sigma(x)$ has full row rank at every $x$.

The following provides a verification of the abstract Lyapunov drift and local oscillation assumptions from the classical diffusion hypotheses.
\begin{proposition}
\label{prop:classverify}
Under Assumption~\ref{ass:SDE}, the Markov semigroup $(P_t)$ satisfies
Assumptions~\ref{ass:L} and \ref{ass:LOC}.
\end{proposition}

\begin{proof}
The generator of \eqref{eq:classicalSDE} is (\cite{StroockVaradhan1979})
\[
 \mathcal Lf(x)
 =
 b(x)\cdot Df(x)
 +
 \frac12
 \tr\left(
  \sigma(x)\sigma(x)^{\mathsf T}D^2f(x)
 \right).
\]
For $V(x)=\abs{x}^2$, condition \eqref{eq:classicalLyap} is precisely
\begin{equation}\label{eq:classical-generator-Lyap}
 \mathcal LV\leq-cV+C_0.
\end{equation}
A standard localisation of It\^o's formula applied to
$e^{ct}V(X_t^x)$ gives
\begin{equation}\label{eq:classicalVbound-detailed}
 P_tV(x)
 \leq
 e^{-ct}V(x)
 +
 \frac{C_0}{c}(1-e^{-ct}),
 \qquad t\geq0.
\end{equation}
Fix $a>0$ and put
\[
 \gamma:=e^{-ca},
 \qquad
 K:=\frac{C_0}{c}(1-e^{-ca}).
\]
Then
\[
 P_aV\leq\gamma V+K.
\]
Moreover, for $0\leq r\leq a$,
\[
 P_rV(x)
 \leq
 V(x)+\frac{C_0}{c}(1-e^{-ca}).
\]
Thus Assumption~\ref{ass:L}, including its short-time part, holds.

We next verify the local oscillation contraction. Choose
\(
 R>\frac{2K}{1-\gamma}
\)
and recall that
\[
 C_R=\{x:V(x)\leq R\}
 =\overline B_{\sqrt R}(0).
\]
Choose a bounded connected open set $O_R$ such that
\[
 C_R\Subset O_R.
\]

On $O_R$, the drift is bounded, and
\eqref{eq:classicalell} gives constants
$0<\lambda_R\leq\Lambda_R<\infty$ such that
\[
 \lambda_RI_d
 \leq
 \sigma(x)\sigma(x)^{\mathsf T}
 \leq
 \Lambda_RI_d,
 \qquad x\in O_R.
\]
Hence
\begin{equation}\label{eq:linearPDE-detailed}
 \partial_tu
 =
 b(x)\cdot Du
 +
 \frac12
 \tr\left(
  \sigma(x)\sigma(x)^{\mathsf T}D^2u
 \right)
\end{equation}
is uniformly parabolic on every space--time cylinder compactly
contained in $(0,\infty)\times O_R$.

Choose $s_*\in(0,a)$ and $x_*\in C_R$. Since $C_R$ is compact and
$O_R$ is connected, a finite family of overlapping spatial balls,
combined with finitely many increasing time levels between $s_*$ and
$a$, gives a parabolic Harnack chain connecting $(s_*,x_*)$ to
$\{a\}\times C_R$. Let $\mathfrak Q_R$ be a bounded positive-time
space--time set containing the cylinders in this chain and satisfying
\[
 \{a\}\times C_R
 \subset
 \mathfrak Q_R
 \Subset
 (0,\infty)\times O_R.
\]

Let $\varphi\in B_V(\R^d)$ and set
\[
 u(t,x):=P_t\varphi(x).
\]
Since
\[
 \abs{\varphi(x)}
 \leq
 \|\varphi\|_V(1+V(x)),
\]
estimate \eqref{eq:classicalVbound-detailed} gives
\[
 \abs{P_t\varphi(x)}
 \leq
 \|\varphi\|_V\bigl(1+P_tV(x)\bigr).
\]
Thus $u$ is bounded on $\mathfrak Q_R$.

For completeness, the passage from bounded to $V$-growing initial
data may be justified by truncation. The global Lipschitz and
linear-growth assumptions imply that, for every $p\geq2$ and every
$T<\infty$,
\[
 \sup_{\substack{0\leq t\leq T\\x\in K}}
 E\bigl[\abs{X_t^x}^p\bigr]
 <\infty
\]
whenever $K\subset\R^d$ is compact. Taking some $p>2$ gives the
uniform integrability needed to pass from bounded truncations of
$\varphi$ to $\varphi$. Consequently, $u$ is a locally bounded
viscosity solution of \eqref{eq:linearPDE-detailed} on the
positive-time region containing $\mathfrak Q_R$.

Define
\[
 m
 :=
 \inf_{(t,x)\in\mathfrak Q_R}u(t,x),
 \qquad
 M
 :=
 \sup_{(t,x)\in\mathfrak Q_R}u(t,x),
\]
and put
\[
 v:=u-m,
 \qquad
 w:=M-u.
\]
Both $v$ and $w$ are nonnegative on $\mathfrak Q_R$. Since
\eqref{eq:linearPDE-detailed} is linear and annihilates constants,
both functions solve the same equation. Iterating the local parabolic
Harnack inequality along the finite chain
 \cite{Krylov1987,Lieberman1996,Wang2011}, we obtain a constant
$C_{H,R}\geq1$, independent of $\varphi$, such that
\[
 v(s_*,x_*)
 \leq
 C_{H,R}\inf_{x\in C_R}v(a,x)
\]
and
\[
 w(s_*,x_*)
 \leq
 C_{H,R}\inf_{x\in C_R}w(a,x).
\]
Lemma~\ref{lem:Hosc} therefore gives
\begin{equation}\label{eq:classical-local-osc}
 \osc_{C_R}(P_a\varphi)
 \leq
 \vartheta_R
 \osc_{\mathfrak Q_R}(P_\cdot\varphi),
 \qquad
 \vartheta_R
 =
 1-\frac{1}{C_{H,R}}
 <1.
\end{equation}
This gives the oscillation estimate required in
Assumption~\ref{ass:LOC}.

Finally, since $\mathfrak Q_R$ is a bounded positive-time set and its
spatial projection is contained in
$O_R\subset\overline B_{\sqrt{\widehat R}}(0)=C_{\widehat R}$ for some
$\widehat R<\infty$, there exists $T_R\geq a$ such that
\[
 \{a\}\times C_R
 \subset
 \mathfrak Q_R
 \subset
 (0,T_R]\times C_{\widehat R}.
\]
Thus Assumption~\ref{ass:LOC}, including its space--time
localisation requirement, is verified. Together with the Lyapunov
estimates established above, this proves that
Assumptions~\ref{ass:L} and \ref{ass:LOC} hold.
\end{proof}
 Thus we obtain the following classical invariant measure and weighted mixing result. 
\begin{theorem}
\label{thm:classical}
Under Assumption~\ref{ass:SDE}, there is a unique invariant probability
measure $\pi$ satisfying
\[
 \int_{\R^d}V\,d\pi<\infty.
\]
Moreover, there exist $C,\rho>0$ such that, for every
$\varphi\in B_b(\R^d)$,
\begin{equation}\label{eq:classicalmix}
 \abs{
  P_t\varphi(x)-\int_{\R^d}\varphi\,d\pi
 }
 \leq
 Ce^{-\rho t}(1+\abs{x}^2)\osc(\varphi),
 \qquad t\geq0.
\end{equation}
\end{theorem}

\begin{proof}
This conclusion is also a direct consequence of the classical Harris
theorem. Indeed, \eqref{eq:classicalVbound-detailed} is a geometric
Foster--Lyapunov condition, while local uniform ellipticity implies
that compact Lyapunov sublevel sets are small. Weighted exponential
ergodicity therefore follows from
\cite{HairerMattingly2011,MeynTweedie2009}.

We briefly explain how the conclusion follows within the present
abstract framework. Proposition~\ref{prop:classverify} and
Corollary~\ref{thm:invariant} give an invariant expectation
$\Et\in\mathfrak E_V$ and
\[
 P_t\varphi(x)\longrightarrow\Et[\varphi],
 \qquad
 \varphi\in B_b(\R^d).
\]
Since every $P_t$ is linear, its pointwise limit $\Et$ is linear on
$B_b(\R^d)$.

The restriction of $\Et$ to $C_0(\R^d)$ is a positive bounded linear
functional. By the Riesz--Markov theorem, it is represented there by
a finite Radon measure $\pi$. The finite Lyapunov moment
$\Et[V]<\infty$ makes this representation tight. More precisely, if
$\chi_R\in C_c(\R^d)$ satisfies
\[
 0\leq\chi_R\leq1,
 \qquad
 \chi_R=1\quad\hbox{on }B_R(0),
\]
then
\[
 0\leq1-\chi_R\leq\frac{V}{R^2},
\]
and hence
\[
 \Et[1-\chi_R]
 \leq
 \frac{\Et[V]}{R^2}.
\]
It follows that $\pi(\R^d)=1$ and
\[
 \int_{\R^d}V\,d\pi\leq\Et[V]<\infty.
\]
The same cutoff argument extends the representation from
$C_0(\R^d)$ to $C_b(\R^d)$.

For every $r>0$, local uniform ellipticity and localisation imply the
positive-time strong Feller property
\[
 P_rB_b(\R^d)\subset C_b(\R^d).
\]
For completeness, this can be seen by applying
Lemma~\ref{lem:BorelDataSmoothing} below to the singleton family
consisting of the law of the ordinary diffusion. Indeed, the strong
Markov property gives the required stopped dynamic programming
identity. On every bounded space--time cylinder, the drift is bounded,
while Assumption~\ref{ass:SDE} gives boundedness and uniform
ellipticity of the covariance matrix. Consequently, for every
$\varphi\in B_b(\R^d)$, the map
\[
 (t,x)\longmapsto P_t\varphi(x)
\]
is locally H\"older continuous on $(0,\infty)\times\R^d$, which proves
the displayed inclusion.

We next identify the invariant expectation with integration against
$\pi$. For $f\in C_b(\R^d)$ and $t\geq0$, invariance of $\Et$ gives
\[
 \int_{\R^d}P_tf\,d\pi
 =
 \Et[P_tf]
 =
 \Et[f]
 =
 \int_{\R^d}f\,d\pi.
\]
Since bounded continuous functions determine finite Borel measures on
$\R^d$, it follows that $\pi P_t=\pi$; thus $\pi$ is invariant.

Now let $\varphi\in B_b(\R^d)$ and fix $r>0$. By the strong Feller
property, $P_r\varphi\in C_b(\R^d)$. Hence, using successively the
invariance of $\Et$, its representation on $C_b(\R^d)$, and the
invariance of $\pi$, we obtain
\[
 \Et[\varphi]
 =
 \Et[P_r\varphi]
 =
 \int_{\R^d}P_r\varphi\,d\pi
 =
 \int_{\R^d}\varphi\,d\pi.
\]
Thus $\Et$ is integration against $\pi$ on all of $B_b(\R^d)$.

If $\nu$ is another invariant probability measure with
$\int V\,d\nu<\infty$, then
\[
 \varphi\longmapsto\int_{\R^d}\varphi\,d\nu
\]
is an invariant element of $\mathfrak E_V$. The uniqueness part of the
abstract contraction theorem therefore gives $\nu=\pi$. Finally,
\eqref{eq:classicalmix} follows from
\eqref{eq:continuousrate} with $V(x)=\abs{x}^2$.

\end{proof}

\begin{remark}
The Harnack step uses only local boundedness of the drift. Thus the
frequently imposed condition $b\in C_b^2$ is unnecessarily strong
here. For example, a globally Lipschitz dissipative drift may be
unbounded and have linear growth.

In the classical linear setting, $u-m$ and $M-u$ solve the same
parabolic equation. For a nonlinear upper semigroup this is no longer
true: \(u-m\) solves the upper equation, whereas \(M-u\) solves the
conjugate lower equation. This is the reason for retaining the
classical argument before treating the \(G\)-SDE cases.
\end{remark}

\section{Borel-data smoothing under a stopped DPP}
\label{sec:BorelDataSmoothing}

The following regularisation result is not specific to a single
diffusion. It applies to families of uniformly nondegenerate It\^o
laws satisfying a stopped dynamic programming principle and will
be used in both of the \(G\)-SDE settings considered below.

For \(z_0=(t_0,x_0)\) and \(r>0\), write
\[
 \mathcal C_r(z_0)
 :=
 (t_0-r^2,t_0]\times B_r(x_0),
\]
and equip space--time with the parabolic distance
\[
 d_{\mathrm p}\bigl((t,x),(s,y)\bigr)
 :=
 |x-y|+|t-s|^{1/2}.
\]

\begin{lemma}
\label{lem:BorelDataSmoothing}
Let \(D\subset (0,\infty)\times\mathbb R^d\) be open, and let
\(v:D\to\mathbb R\) be bounded and universally measurable. For
\(z_0=(t_0,x_0)\) and \(r>0\), suppose that whenever
\(
\overline{\mathcal C_{4r}(z_0)}\subset D
\)
and \(z=(t,x)\in \mathcal C_{4r}(z_0)\), there is a nonempty family
\(\mathfrak P(t,x)\) of probability laws such that, under every
\(P\in\mathfrak P(t,x)\),
\[
X_s^x
=
x+\int_0^s\beta_\ell^P\,d\ell+M_s^P,
\]
up to the exit time of the backward space--time process
\(
Z_s^{t,x}:=(t-s,X_s^x)
\)
from \(\mathcal C_{4r}(z_0)\). Here \(M^P\) is a continuous local martingale and
\[
|\beta_\ell^P|\le B,
\qquad
\lambda I_d
\le
\frac{d\langle M^P\rangle_\ell}{d\ell}
\le
\Lambda I_d,
\qquad
d\ell\otimes P\text{-a.e.},
\]
for constants \(0<\lambda\le\Lambda<\infty\), uniformly over the
admissible laws.

Assume, moreover, the following stopped dynamic programming property:
if \(C\) is a cylinder satisfying
\[
\overline C\subset \mathcal C_{4r}(z_0),
\qquad
(t,x)\in C,
\]
and if \(\theta\) is a stopping time bounded by the exit time
\(\tau_C^{t,x}\) of \(Z^{t,x}\) from \(C\), then
\begin{equation}
\label{eq:stopped-DPP}
v(t,x)
=
\sup_{P\in\mathfrak P(t,x)}
E^P\!\left[
v(t-\theta,X_\theta^x)
\right].
\end{equation}

Then there exist \(\alpha\in(0,1)\) and \(C<\infty\), depending only on
\[
d,\qquad \lambda,\qquad \Lambda,\qquad rB,
\]
and on the fixed relative geometry of the cylinders, such that
\begin{equation}
\label{eq:BorelHolder}
|v(z)-v(z')|
\le
C\,\operatorname{osc}_{\mathcal C_{4r}(z_0)}v
\left(
\frac{d_p(z,z')}{r}
\right)^\alpha,
\qquad
z,z'\in \mathcal{C}_r(z_0),
\end{equation}
where
\[
d_p\bigl((t,x),(s,y)\bigr)
:=
|x-y|+|t-s|^{1/2}.
\]
In particular, the pointwise-defined function \(v\) is locally
parabolically H\"older continuous in \(D\).
\end{lemma}

\begin{proof}
The proof is the probabilistic Krylov--Safonov oscillation argument.
We first record the uniform hitting estimate that will be used.

Let \(\zeta=(s_0,y_0)\), let \(0<\rho\le r\), and suppose that
\[
\overline{\mathcal{C}_{4\rho}(\zeta)}
\subset \mathcal C_{4r}(z_0).
\]
Introduce the earlier intermediate cylinder
\[
A_\rho(\zeta)
:=
(s_0-6\rho^2,s_0-5\rho^2)
\times B_\rho(y_0).
\]
Thus
\[
A_\rho(\zeta)\Subset \mathcal C_{3\rho}(\zeta),
\]
and \(A_\rho(\zeta)\) lies a parabolic distance of order \(\rho\)
strictly before every point of \(\mathcal C_\rho(\zeta)\).

The Krylov--Safonov hitting estimate implies that there exists
\(p_0\in(0,1)\), depending only on
\[
d,\qquad \lambda,\qquad \Lambda,\qquad rB,
\]
and on the above fixed relative geometry, with the following property.
If \(\Gamma\subset A_\rho(\zeta)\) is closed and
\begin{equation}
\label{eq:GammaMeasure}
|\Gamma|\ge \frac14 |A_\rho(\zeta)|,
\end{equation}
then
\begin{equation}
\label{eq:uniform-hitting}
\inf_{\substack{(t,x)\in \mathcal{C}_\rho(\zeta)\\
                 P\in\mathfrak P(t,x)}}
P\!\left(
\sigma_\Gamma<\tau_{3\rho}
\right)
\ge p_0.
\end{equation}
Here
\[
\sigma_\Gamma
:=
\inf\{s\ge0:Z_s^{t,x}\in\Gamma\},
\qquad
\tau_{3\rho}
:=
\inf\{s\ge0:Z_s^{t,x}\notin \mathcal C_{3\rho}(\zeta)\}.
\]
The constant is uniform over \(P\) because the probabilistic
Krylov--Safonov estimate depends only on the bounds for the differential
characteristics. Indeed, after the parabolic scaling
\[
s=\rho^2\widetilde s,
\qquad
X_s-y_0=\rho\widetilde X_{\widetilde s},
\]
the ellipticity constants remain \(\lambda,\Lambda\), whereas the
rescaled drift is bounded by \(\rho B\le rB\). See, for example, the
Krylov--Safonov hitting theorem for uniformly nondegenerate It\^o
processes (\cite{Krylov2021Diffusions, KrylovSafonov1979}).

We now prove oscillation decay. Set
\[
M_{4\rho}
:=
\sup_{\mathcal C_{4\rho}(\zeta)}v,
\qquad
m_{4\rho}
:=
\inf_{\mathcal C_{4\rho}(\zeta)}v,
\qquad
\mu_{4\rho}
:=
\frac{M_{4\rho}+m_{4\rho}}2.
\]
Since
\[
A_\rho(\zeta)
=
\bigl(A_\rho(\zeta)\cap\{v\le\mu_{4\rho}\}\bigr)
\cup
\bigl(A_\rho(\zeta)\cap\{v\ge\mu_{4\rho}\}\bigr),
\]
at least one of these two level sets has Lebesgue measure at least
\(\frac12|A_\rho(\zeta)|\). The level sets are Lebesgue measurable
because \(v\) is universally measurable. By inner regularity, the
larger level set contains a compact subset \(\Gamma\) satisfying
\eqref{eq:GammaMeasure}.

Suppose first that
\[
\Gamma\subset\{v\le\mu_{4\rho}\}.
\]
Fix \((t,x)\in \mathcal{C}_\rho(\zeta)\) and set
\[
\theta:=\sigma_\Gamma\wedge\tau_{3\rho}.
\]
Since
\[
\overline{\mathcal C_{3\rho}(\zeta)}
\subset \mathcal C_{4\rho}(\zeta),
\]
the stopped point \(Z_\theta^{t,x}\) belongs to
\(\mathcal{C}_{4\rho}(\zeta)\). Hence, under every
\(P\in\mathfrak P(t,x)\),
\[
v(t-\theta,X_\theta^x)
\le
\mu_{4\rho}
  I_{\{\sigma_\Gamma<\tau_{3\rho}\}}
+
M_{4\rho}
  I_{\{\sigma_\Gamma\ge\tau_{3\rho}\}}.
\]
Using \eqref{eq:uniform-hitting}, we obtain
\[
\begin{aligned}
E^P\!\left[v(t-\theta,X_\theta^x)\right]
&\le
p_0\mu_{4\rho}+(1-p_0)M_{4\rho} \\
&=
M_{4\rho}
-\frac{p_0}{2}
 \bigl(M_{4\rho}-m_{4\rho}\bigr).
\end{aligned}
\]
Taking the supremum over \(P\) in the stopped DPP gives
\[
v(t,x)
\le
M_{4\rho}
-\frac{p_0}{2}
 \operatorname{osc}_{\mathcal C_{4\rho}(\zeta)}v.
\]
Therefore,
\begin{equation}
\label{eq:upper-improvement}
\sup_{\mathcal{C}_\rho(\zeta)}v
\le
M_{4\rho}
-\frac{p_0}{2}
 \operatorname{osc}_{\mathcal C_{4\rho}(\zeta)}v.
\end{equation}

Suppose instead that
\[
\Gamma\subset\{v\ge\mu_{4\rho}\}.
\]
For every \(P\in\mathfrak P(t,x)\),
\[
v(t-\theta,X_\theta^x)
\ge
\mu_{4\rho}
  I_{\{\sigma_\Gamma<\tau_{3\rho}\}}
+
m_{4\rho}
  I_{\{\sigma_\Gamma\ge\tau_{3\rho}\}}.
\]
Consequently,
\[
\begin{aligned}
E^P\!\left[v(t-\theta,X_\theta^x)\right]
&\ge
p_0\mu_{4\rho}+(1-p_0)m_{4\rho}\\
&=
m_{4\rho}
+\frac{p_0}{2}
 \bigl(M_{4\rho}-m_{4\rho}\bigr).
\end{aligned}
\]
Since the stopped DPP is a supremum over \(P\), this lower bound,
which is valid for every admissible \(P\), gives
\[
v(t,x)
\ge
m_{4\rho}
+\frac{p_0}{2}
 \operatorname{osc}_{\mathcal C_{4\rho}(\zeta)}v.
\]
Hence
\begin{equation}
\label{eq:lower-improvement}
\inf_{\mathcal{C}_\rho(\zeta)}v
\ge
m_{4\rho}
+\frac{p_0}{2}
 \operatorname{osc}_{\mathcal C_{4\rho}(\zeta)}v.
\end{equation}

In either case, \eqref{eq:upper-improvement} or
\eqref{eq:lower-improvement} yields
\begin{equation}
\label{eq:oscillation-decay}
\operatorname{osc}_{\mathcal C_\rho(\zeta)}v
\le
\vartheta\,
\operatorname{osc}_{\mathcal C_{4\rho}(\zeta)}v,
\qquad
\vartheta:=1-\frac{p_0}{2}\in(0,1).
\end{equation}

Iterating \eqref{eq:oscillation-decay} on cylinders whose radii decrease
by factors of \(4\) gives
\[
\operatorname{osc}_{\mathcal C_\rho(\zeta)}v
\le
C
\left(\frac{\rho}{r}\right)^\alpha
\operatorname{osc}_{\mathcal C_{4r}(z_0)}v,
\qquad
0<\rho\le r,
\]
where one may take
\[
\alpha
=
\min\left\{
\frac {1}{2},\frac{-\log \vartheta}{\log 4}
\right\}.
\]
A standard covering argument in the parabolic metric, using the buffer
between \(\mathcal C_r(z_0)\) and \(\mathcal C_{4r}(z_0)\), now gives
\eqref{eq:BorelHolder} for arbitrary \(z,z'\in \mathcal C_r(z_0)\).

Since the estimate applies directly to the pointwise values of \(v\),
it proves that \(v\) itself is locally H\"older continuous; no
modification or semicontinuous representative is required.
\end{proof}
\section{Multiplicative \texorpdfstring{$G$}{G}-SDEs}
\label{sec:multG}

Let $B$ be an $m$-dimensional $G$-Brownian motion. We use the
normalisation
\begin{equation}\label{eq:Gdef}
 G(A)=\frac12\sup_{Q\in\Q}\tr(QA),
 \qquad A\in\mathbb S^m,
\end{equation}
where $\Q\subset\mathbb S_+^m$ is nonempty, compact and convex. We
consider the multiplicative $G$-SDE
\begin{equation}\label{eq:multGSDE}
 dX_t^x=b(X_t^x)\,dt+\sigma(X_t^x)\,dB_t,
 \qquad X_0^x=x.
\end{equation}
Throughout this section, $b:\R^d\to\R^d$ and
$\sigma:\R^d\to\R^{d\times m}$ are globally Lipschitz and have at
most linear growth. These assumptions are repeated in
Assumption~\ref{ass:multG} for ease of reference.

Let
\[
 \Omega:=C_0(\mathbb R_+;\mathbb R^m),
 \qquad
 B_t(\omega):=\omega_t,
\]
be the canonical path space and its coordinate process. For \(T>0\), let
\[
\Omega_T:=C_0([0,T];\mathbb R^m), 
\] and
let \(L_{ip}(\Omega_T)\) denote the space of bounded Lipschitz
cylinder functions of the form
\[
 \xi
 =
 \psi(B_{t_1},\ldots,B_{t_n}),
 \qquad
 0\leq t_1<\cdots<t_n\leq T,
\]
where
\[
 \psi\in
 C_{b,\mathrm{Lip}}\bigl((\mathbb R^m)^n\bigr).
\]
The \(G\)-expectation \(\widehat{\mathbb E}\) is initially defined on
\(L_{ip}(\Omega_T)\). We write \(L_G^1(\Omega_T)\) for its completion
under the norm
\[
 \|\xi\|_{L_G^1}
 :=
 \widehat{\mathbb E}[|\xi|],
\]
with functions of zero norm identified. Similarly,
\(L_G^1(\Omega)\) is the completion of
\[
 L_{ip}(\Omega)
 :=
 \bigcup_{T>0}L_{ip}(\Omega_T)
\]
under the same norm.

On the natural $G$-expectation domain, the associated semigroup is
\[
 \mathcal T_t\varphi(x)
 :=\widehat{\mathbb E}\bigl[\varphi(X_t^x)\bigr].
\]
For bounded Lipschitz $\varphi$, the representation theorem for
$G$-expectation identifies this operator with the usual strong-control
value function. More precisely, let
$(\Omega^0,\mathcal F^0,\mathbb F^0,P_0)$ carry a standard
$m$-dimensional Brownian motion $W$, and let $\mathcal A_t$ be the
collection of progressively measurable processes $Q_\cdot$ on
$[0,t]$ taking values in $\Q$. For $Q_\cdot\in\mathcal A_t$, set
\begin{equation}\label{eq:BQdefinition}
 B_s^Q:=\int_0^sQ_r^{1/2}\,dW_r,
 \qquad 0\le s\le t,
\end{equation}
so that
\begin{equation}\label{eq:BQquadraticvariation}
 d\langle B^Q\rangle_s=Q_s\,ds.
\end{equation}
The corresponding state equation is
\begin{equation}\label{eq:multcontrolBQ}
 dX_s^{x,Q}
 =b(X_s^{x,Q})\,ds+\sigma(X_s^{x,Q})\,dB_s^Q,
 \qquad X_0^{x,Q}=x,
\end{equation}
or, equivalently,
\begin{equation}\label{eq:multcontrol}
 dX_s^{x,Q}
 =b(X_s^{x,Q})\,ds
  +\sigma(X_s^{x,Q})Q_s^{1/2}\,dW_s,
 \qquad X_0^{x,Q}=x.
\end{equation}
The representation theorem of Denis--Hu--Peng
\cite{DenisHuPeng2011} states that, for every
\(\xi\in L_G^1(\Omega_t)\),
\begin{equation}\label{eq:Grepresentationcontrol}
 \widehat{\mathbb E}[\xi]
 =
 \sup_{Q_\cdot\in\mathcal A_t}
 E_{P_0}\bigl[\xi(B^Q)\bigr].
\end{equation}
Under the pullback $B\mapsto B^Q$, pathwise uniqueness gives
\begin{equation}\label{eq:pullbacksolution}
 X_s^x(B^Q)=X_s^{x,Q},
 \qquad 0\le s\le t,\quad P_0\text{-a.s.}
\end{equation}
Consequently, whenever
\(\varphi(X_t^x)\in L_G^1(\Omega_t)\),
\begin{equation}\label{eq:agreementGcontrol}
 \mathcal T_t\varphi(x)
 =\sup_{Q_\cdot\in\mathcal A_t}
   E_{P_0}\bigl[\varphi(X_t^{x,Q})\bigr].
\end{equation}
In particular, this holds for
$\varphi\in C_{b,\mathrm{Lip}}(\R^d)$.

For general bounded Borel payoffs, we use instead the canonical-law
extension below. We do not require the strong-control expression in
\eqref{eq:agreementGcontrol} to remain valid outside the natural
$G$-expectation domain.

\subsection{Canonical-law extension and dynamic programming}

Let
\[
 \Omega:=C_0(\mathbb R_+;\mathbb R^m)
\]
be equipped with the topology of locally uniform convergence, and let
\[
 \mathcal F:=\mathcal B(\Omega).
\]
Define the canonical coordinate process by
\[
 B_t(\omega):=\omega(t),
 \qquad t\geq0,\quad \omega\in\Omega,
\]
and let
\[
 \mathcal F_t
 :=
 \sigma(B_s:0\leq s\leq t),
 \qquad t\geq0.
\]
We write
\[
 \mathbb F=(\mathcal F_t)_{t\geq0}
\]
for the raw canonical filtration.
Let $\mathfrak P(\Omega)$ denote the space of probability measures on
$\Omega$ with its weak topology. Following Nutz--van Handel
\cite{NutzVanHandel2013}  and Nutz \cite{Nutz2012, Nutz2013}, let $\langle B\rangle$ be a pathwise,
progressively measurable version of the quadratic variation and let
$\widehat Q_t=d\langle B\rangle_t/dt$ whenever the latter exists.
Define
\begin{equation}\label{eq:canonicalvolatilitylaws}
 \mathcal P_{\Q}
 :=
 \left\{
 P\in\mathfrak P(\Omega):
 \begin{array}{l}
 B\text{ is a continuous local $P$-martingale},\\
 \langle B\rangle\text{ is $P$-a.s. absolutely continuous},\\
 \widehat Q_t\in\Q\quad dt\otimes P\text{-a.e.}
 \end{array}
 \right\}.
\end{equation}
Since $\Q$ is bounded, $B$ is a true martingale under every
$P\in\mathcal P_{\Q}$.

We shall use the following standard universal version of the state
equation. Its role is only to ensure that the payoff entering the
canonical conditional expectation is a measurable path functional.

\begin{lemma}
\label{lem:universal-state-map}
There exists a map
\[
(x,\omega,t)\longmapsto X_t^x(\omega)
\]
from
\[
\mathbb R^d\times\Omega\times\mathbb R_+
\quad\text{to}\quad
\mathbb R^d
\]
which is
\[
\mathcal B(\mathbb R^d)\otimes\operatorname{Prog}(\mathbb F)
\text{-measurable},
\]
such that, for every \(P\in\mathcal P_{\mathcal Q}\) and every
\(x\in\mathbb R^d\),
\begin{equation}
\label{eq:universal-SDE}
X_t^x
=
x+\int_0^t b(X_s^x)\,ds
  +\int_0^t\sigma(X_s^x)\,dB_s,
\qquad t\ge0,
\quad P\text{-a.s.}
\end{equation}
In particular, \(X^x\) has continuous paths \(P\)-a.s. and is
indistinguishable from the unique \(P\)-solution of
\eqref{eq:universal-SDE}.

Moreover, the following conditional flow property holds. 
For every $P\in \mathcal P_{\mathcal Q}$, every finite stopping time $\tau$, and every admissible continuation kernel \(\nu(\omega)\in\mathcal P_{\mathcal Q}\), the conditional flow property holds under the corresponding pasted law
\begin{equation}
\label{eq:conditional-flow}
X_{\tau(\omega)+s}^x
   \bigl(\omega\otimes_\tau\widetilde\omega\bigr)
=
X_s^{X_\tau^x(\omega)}(\widetilde\omega),
\qquad s\ge0,
\end{equation}
for
\[
P(d\omega)\nu(\omega,d\widetilde\omega)\text{-a.e. }
(\omega,\widetilde\omega).
\]
Thus the flow property holds in precisely the almost-sure sense
required by conditioning and pasting.
\end{lemma}

\begin{proof}
For \(n\ge1\), set
\[
\Delta_n:=2^{-n},
\qquad
\eta_n(s):=
\Delta_n\left\lfloor\frac{s}{\Delta_n}\right\rfloor,
\]
and define the continuous-time Euler approximation by
\begin{equation}
\label{eq:Euler-universal}
X_t^{n,x}
:=
x+\int_0^t b\bigl(X_{\eta_n(s)}^{n,x}\bigr)\,ds
  +\int_0^t
       \sigma\bigl(X_{\eta_n(s)}^{n,x}\bigr)\,dB_s.
\end{equation}
The stochastic integral in \eqref{eq:Euler-universal} is a finite sum
of increments of the canonical process on each bounded time interval.
Consequently,
\[
(x,\omega,t)\longmapsto X_t^{n,x}(\omega)
\]
is jointly Borel measurable, and for every fixed \(x\), the map
\((\omega,t)\mapsto X_t^{n,x}(\omega)\) is progressively measurable. 

Fix \(P\in\mathcal P_{\mathcal Q}\). Let \(X^{x,P}\) denote the unique
\(P\)-solution of
\[
X_t^{x,P}
=
x+\int_0^t b(X_s^{x,P})\,ds
  +\int_0^t\sigma(X_s^{x,P})\,dB_s.
\]
Since
\[
d\langle B\rangle_s=\widehat Q_s^P\,ds,
\qquad
\widehat Q_s^P\in\mathcal Q,
\quad ds\otimes P\text{-a.e.},
\]
compactness of \(\mathcal Q\) gives a constant
\(\overline q<\infty\) such that
\[
0\le \widehat Q_s^P\le\overline q I_m,
\qquad ds\otimes P\text{-a.e.}
\]

The standard BDG and Gronwall estimates imply that, for every
\(p\ge2\), \(T>0\), and compact \(K\subset\mathbb R^d\),
\begin{equation}
\label{eq:universal-moment-bound}
\sup_{\substack{P\in\mathcal P_{\mathcal Q}\\x\in K}}
E^P\left[
\sup_{0\le s\le T}
\bigl(
|X_s^{x,P}|^p+|X_s^{n,x}|^p
\bigr)
\right]
\le C_{p,T,K}.
\end{equation}
Moreover,
\begin{equation}
\label{eq:Euler-increment-bound}
\sup_{\substack{P\in\mathcal P_{\mathcal Q}\\x\in K}}
E^P\left[
\left|
X_s^{n,x}-X_{\eta_n(s)}^{n,x}
\right|^p
\right]
\le
C_{p,T,K}\Delta_n^{p/2},
\qquad 0\le s\le T.
\end{equation}
Global Lipschitz continuity of \(b\) and \(\sigma\), followed by BDG,
Gronwall's inequality, and
\eqref{eq:Euler-increment-bound}, yields
\begin{equation}
\label{eq:Euler-convergence-bound}
\sup_{\substack{P\in\mathcal P_{\mathcal Q}\\x\in K}}
E^P\left[
\sup_{0\le s\le T}
|X_s^{n,x}-X_s^{x,P}|^p
\right]
\le
C_{p,T,K}\Delta_n^{p/2}
=
C_{p,T,K}2^{-np/2}.
\end{equation}

For fixed \(P,x,T\), and \(\varepsilon>0\), Markov's inequality gives
\[
\sum_{n=1}^{\infty}
P\left(
\sup_{0\le s\le T}
|X_s^{n,x}-X_s^{x,P}|>\varepsilon
\right)
\le
C_{p,T,K}\varepsilon^{-p}
\sum_{n=1}^{\infty}2^{-np/2}
<\infty.
\]
Hence the Borel--Cantelli lemma implies
\begin{equation}
\label{eq:Euler-as-convergence}
X^{n,x}\longrightarrow X^{x,P}
\quad\text{uniformly on }[0,T],
\qquad P\text{-a.s.}
\end{equation}
Taking the intersection over \(T\in\mathbb N\), we obtain locally
uniform convergence on \(\mathbb R_+\), \(P\)-a.s.

We now define the universal map. For \(i=1,\ldots,d\), put
\[
L_t^{x,i}(\omega)
:=
\limsup_{n\to\infty}X_t^{n,x,i}(\omega).
\]
Set
\begin{equation}
\label{eq:progressive-limsup-definition}
X_t^x(\omega)
:=
\begin{cases}
\bigl(L_t^{x,1}(\omega),\ldots,L_t^{x,d}(\omega)\bigr),
 &\text{if every }L_t^{x,i}(\omega)\text{ is finite},\\[1mm]
x,
 &\text{otherwise}.
\end{cases}
\end{equation}
This is the crucial change from the earlier construction. The
definition at time \(t\) uses only the values
\(\{X_t^{n,x}\}_{n\ge1}\) at that same time. It does not use the event
that the entire sequence converges on all future time intervals.
Since the progressive sigma-field is closed under pointwise
\(\limsup\), \eqref{eq:progressive-limsup-definition} is
\[
\mathcal B(\mathbb R^d)\otimes\operatorname{Prog}(\mathbb F)
\text{-measurable}.
\]
On the \(P\)-full set in \eqref{eq:Euler-as-convergence}, the
coordinatewise \(\limsup\) is the locally uniform limit. Therefore,
for every \(P\in\mathcal P_{\mathcal Q}\) and every
\(x\in\mathbb R^d\),
\[
X^x=X^{x,P},
\qquad P\text{-a.s.}
\]
Thus \(X^x\) satisfies \eqref{eq:universal-SDE} and has continuous
paths \(P\)-a.s.
It remains to verify the flow property. Let \(\tau\) be a finite
stopping time and let \(\nu\) be an admissible continuation kernel.
By the conditioning and pasting stability of
\(\mathcal P_{\mathcal Q}\), the concatenated canonical process is a
continuous local martingale whose quadratic-variation density remains
\(\mathcal Q\)-valued.
For fixed past path \(\omega\), define
\[
Y_s(\widetilde\omega)
:=
X_{\tau(\omega)+s}^x
 \bigl(\omega\otimes_\tau\widetilde\omega\bigr).
\]
Conditioning \eqref{eq:universal-SDE} at \(\tau\) shows that, for
\(P\)-a.e. \(\omega\), the process \(Y\) satisfies, under
\(\nu(\omega)\),
\[
Y_s
=
X_\tau^x(\omega)
+\int_0^s b(Y_r)\,dr
+\int_0^s\sigma(Y_r)\,d\widetilde B_r.
\]
On the other hand, by the already established universal property,
\[
Z_s(\widetilde\omega)
:=
X_s^{X_\tau^x(\omega)}(\widetilde\omega)
\]
satisfies the same SDE under \(\nu(\omega)\), with the same initial
condition and the same canonical driving process. Global Lipschitz
continuity gives pathwise uniqueness. Hence
\[
Y_s=Z_s,\qquad s\ge0,
\quad \nu(\omega)\text{-a.s.}
\]
for \(P\)-a.e. \(\omega\). This is exactly
\eqref{eq:conditional-flow}.
\end{proof}
\begin{definition} Let \(S\) be a Polish space. A subset \(A\subset S\) is called
\emph{analytic} if there exist a Polish space \(Y\) and a Borel set
\(B\subset S\times Y\) such that
\[
 A=\operatorname{proj}_S(B).
\]
A function \(f:S\to[-\infty,\infty]\) is called
\emph{upper semianalytic} if
\[
 \{x\in S:f(x)>c\}
\]
is analytic for every \(c\in\mathbb R\). Every Borel function is
upper semianalytic, and every upper semianalytic function is
universally measurable.
\end{definition}
The upper semianalytic class is used because it is preserved by
integration against a variable probability measure and by taking
suprema over correspondences with analytic graphs. Borel
measurability, in contrast, is not generally preserved by such
uncountable optimisation.
\begin{lemma}
\label{lem:joint-value-measurability}
Let
\[
 \varphi:\mathbb R^d\longrightarrow\mathbb R
\]
be bounded and upper semianalytic. Then the maps
\[
 (t,x,\omega)\longmapsto \varphi(X_t^x(\omega))
\]
and
\[
 (t,x,P)\longmapsto E^P[\varphi(X_t^x)]
\]
are upper semianalytic on
\[
 \mathbb R_+\times\mathbb R^d\times\Omega
 \quad\text{and}\quad
 \mathbb R_+\times\mathbb R^d\times\mathfrak P(\Omega),
\]
respectively.
Moreover, the set
\(\mathcal P_{\mathcal Q}\subset\mathfrak P(\Omega)\)
defined in \eqref{eq:canonicalvolatilitylaws} is Borel. Consequently,
\[
 (t,x)\longmapsto
 T_t\varphi(x)
 :=
 \sup_{P\in\mathcal P_{\mathcal Q}}
 E^P[\varphi(X_t^x)]
\]
is upper semianalytic on
\(\mathbb R_+\times\mathbb R^d\), and hence universally measurable.
\end{lemma}
\begin{proof}
Set
\[
 Z:=\mathbb R_+\times\mathbb R^d.
\]
Since \(\Omega\), \(Z\), and \(\mathfrak P(\Omega)\) are Polish spaces,
all the products appearing below are standard Borel spaces.
\medskip
\noindent
\emph{Step 1: measurability of the payoff.}
By Lemma~\ref{lem:universal-state-map}, the map
\[
 (x,\omega,t)\longmapsto X_t^x(\omega)
\]
is
\[
 \mathcal B(\mathbb R^d)
 \otimes\operatorname{Prog}(\mathbb F)
 \text{-measurable}.
\]
Since
\[
 \operatorname{Prog}(\mathbb F)
 \subset
 \mathcal B(\mathbb R_+)\otimes\mathcal F,
\]
the map
\[
 F:Z\times\Omega\longrightarrow\mathbb R^d,
 \qquad
 F(t,x,\omega):=X_t^x(\omega),
\]
is Borel measurable.
For every \(c\in\mathbb R\),
\[
 \bigl\{(t,x,\omega):
       \varphi(X_t^x(\omega))>c\bigr\}
 =
 F^{-1}\bigl(\{y:\varphi(y)>c\}\bigr).
\]
Because \(\varphi\) is upper semianalytic,
\(\{y:\varphi(y)>c\}\) is analytic. The inverse image of an
analytic set under a Borel map between standard Borel spaces is
analytic. Hence
\[
 (t,x,\omega)\longmapsto\varphi(X_t^x(\omega))
\]
is upper semianalytic.
\medskip
\noindent
\emph{Step 2: integration with respect to a variable probability
measure.}
Define a probability kernel from
\(Z\times\mathfrak P(\Omega)\) to \(\Omega\) by
\[
 K((t,x,P),A):=P(A),
 \qquad A\in\mathcal F.
\]
This is a Borel kernel: for every \(A\in\mathcal F\), the evaluation
map
\[
 P\longmapsto P(A)
\]
is Borel on \(\mathfrak P(\Omega)\).
A Borel kernel maps bounded upper semianalytic integrands to upper
semianalytic parameterised integrals. Applying this fact to
\[
 g(t,x,\omega):=\varphi(X_t^x(\omega))
\]
gives that
\[
 H(t,x,P)
 :=
 \int_\Omega g(t,x,\omega)\,P(d\omega)
 =
 E^P[\varphi(X_t^x)]
\]
is upper semianalytic on
\(Z\times\mathfrak P(\Omega)\).
\medskip
\noindent
\emph{Step 3: the set \(\mathcal P_{\mathcal Q}\) is Borel.}
Let
\[
 \mathfrak M_a
 :=
 \left\{
 P\in\mathfrak P(\Omega):
 \begin{array}{l}
 B\text{ is a continuous local \(P\)-martingale},\\
 \langle B\rangle\text{ is \(P\)-a.s. absolutely continuous}
 \end{array}
 \right\}.
\]
By \cite[Lemma~4.4]{NutzVanHandel2013},
\(\mathfrak M_a\) is a Borel subset of
\(\mathfrak P(\Omega)\).
We use the pathwise, progressively measurable quadratic variation
\(\langle B\rangle\) and the progressively measurable density
constructed in \cite[(4.2)]{NutzVanHandel2013}. More precisely, for
\(s>0\), define componentwise
\[
 \overline Q_s(\omega)
 :=
 \limsup_{n\to\infty}
 n\left(
   \langle B\rangle_s(\omega)
   -
   \langle B\rangle_{(s-1/n)\vee0}(\omega)
 \right).
\]
Choose a fixed matrix
\(Q_\ast\in\mathbb S^m\setminus\mathcal Q\), and set
\[
 \widehat Q_s(\omega)
 :=
 \begin{cases}
  \overline Q_s(\omega),
    &\text{if all entries of \(\overline Q_s(\omega)\) are finite},\\
  Q_\ast,
    &\text{otherwise},
 \end{cases}
 \qquad
 \widehat Q_0:=Q_\ast.
\]
Then \(\widehat Q\) is defined on all of
\(\mathbb R_+\times\Omega\) and is progressively measurable.
For every \(P\in\mathfrak M_a\),
\[
 \widehat Q_s
 =
 \frac{d\langle B\rangle_s}{ds}
 \qquad ds\otimes P\text{-a.e.}
\]
Define
\[
 A_{\mathcal Q}
 :=
 \bigcap_{N=1}^{\infty}
 \left\{
 \omega\in\Omega:
 \int_0^N
 \mathbf 1_{\{\widehat Q_s(\omega)\notin\mathcal Q\}}\,ds
 =0
 \right\}.
\]
Since \(\widehat Q\) is progressively measurable and
\(\mathcal Q\) is Borel, the map
\[
 (s,\omega)\longmapsto
 \mathbf 1_{\{\widehat Q_s(\omega)\notin\mathcal Q\}}
\]
is Borel measurable. Consequently, for every \(N\ge1\),
\[
 \omega\longmapsto
 \int_0^N
 \mathbf 1_{\{\widehat Q_s(\omega)\notin\mathcal Q\}}\,ds
\]
is Borel measurable. It follows that \(A_{\mathcal Q}\) is a Borel
subset of \(\Omega\).
If \(P\in\mathfrak M_a\), then Tonelli's theorem gives
\[
\begin{aligned}
 P(A_{\mathcal Q})=1
 &\quad\Longleftrightarrow\quad
 \int_0^N
 \mathbf 1_{\{\widehat Q_s(\omega)\notin\mathcal Q\}}\,ds
 =0
 \quad P\text{-a.s. for every }N
\\
 &\quad\Longleftrightarrow\quad
 \int_\Omega\int_0^N
 \mathbf 1_{\{\widehat Q_s(\omega)\notin\mathcal Q\}}
 \,ds\,P(d\omega)
 =0
 \quad\text{for every }N
\\
 &\quad\Longleftrightarrow\quad
 \widehat Q_s(\omega)\in\mathcal Q
 \quad ds\otimes P\text{-a.e.}
\end{aligned}
\]
Thus
\[
 \mathcal P_{\mathcal Q}
 =
 \left\{
 P\in\mathfrak M_a:
 P(A_{\mathcal Q})=1
 \right\}.
\]
Because \(A_{\mathcal Q}\) is Borel, the evaluation map
\[
 P\longmapsto P(A_{\mathcal Q})
\]
is Borel on \(\mathfrak P(\Omega)\). Therefore
\(\mathcal P_{\mathcal Q}\) is Borel.
\medskip
\noindent
\emph{Step 4: optimisation over \(\mathcal P_{\mathcal Q}\).}
Regard \(\mathcal P_{\mathcal Q}\) as the constant correspondence
\[
 \Gamma:Z\longrightarrow 2^{\mathfrak P(\Omega)},
 \qquad
 \Gamma(t,x):=\mathcal P_{\mathcal Q}.
\]
Its graph is
\[
 \operatorname{Gr}(\Gamma)
 =
 Z\times\mathcal P_{\mathcal Q}.
\]
By Step~3, this graph is Borel and hence analytic.
For \(c\in\mathbb R\), we have
\[
\begin{aligned}
 \bigl\{(t,x):T_t\varphi(x)>c\bigr\}
 &=
 \operatorname{proj}_Z
 \left(
  \left\{
   (t,x,P)\in\operatorname{Gr}(\Gamma):
   H(t,x,P)>c
  \right\}
 \right).
\end{aligned}
\]
Since \(H\) is upper semianalytic, the set
\[
 \{(t,x,P):H(t,x,P)>c\}
\]
is analytic. Its intersection with the Borel set
\(\operatorname{Gr}(\Gamma)\) is analytic, and the projection of an
analytic set is analytic. Therefore
\[
 \{(t,x):T_t\varphi(x)>c\}
\]
is analytic for every \(c\in\mathbb R\). This proves that
\((t,x)\mapsto T_t\varphi(x)\) is upper semianalytic and hence
universally measurable.
\end{proof}
\begin{remark}
\label{rem:stopped-value-measurability}
Let \(t\ge0\) and let \(\theta\le t\) be a bounded stopping time.
Since \(X^x\) is progressively measurable, the map
\[
 \omega\longmapsto
 \bigl(t-\theta(\omega),X_{\theta(\omega)}^x(\omega)\bigr)
\]
is measurable. Lemma~\ref{lem:joint-value-measurability} therefore
implies that
\[
 \omega\longmapsto
 T_{t-\theta(\omega)}
 \varphi\bigl(X_{\theta(\omega)}^x(\omega)\bigr)
\]
is universally measurable; more precisely, it is measurable with
respect to the universal completion of
\(\mathcal F_\theta\). Thus this random variable is well defined under
every \(P\in\mathcal P_{\mathcal Q}\).
\end{remark}
For every bounded upper semianalytic function
$\varphi:\R^d\to\R$, define
\begin{equation}\label{eq:multT}
 T_t\varphi(x)
 :=\sup_{P\in\mathcal P_{\Q}}
   E^P\bigl[\varphi(X_t^x)\bigr].
\end{equation}
The canonical representation of $G$-expectation gives
\[
 \sup_{P\in\mathcal P_{\Q}}E^P[\xi]
 =
 \widehat{\mathbb E}[\xi],
 \qquad
 \xi\in L_G^1(\Omega_t);
\]
see \cite[Section~3]{NutzVanHandel2013}. Consequently, whenever
$\varphi(X_t^x)\in L_G^1(\Omega_t)$,
\[
 T_t\varphi(x)
 =
 \sup_{P\in\mathcal P_{\Q}}
 E^P\bigl[\varphi(X_t^x)\bigr]
 =
 \widehat{\mathbb E}\bigl[\varphi(X_t^x)\bigr]
 =
 \mathcal T_t\varphi(x).
\]
Thus the canonical-law semigroup $T_t$ and the original
$G$-expectation semigroup $\mathcal T_t$ agree on their common domain.
In particular, they agree for
$\varphi\in C_{b,\mathrm{Lip}}(\R^d)$, and this common value also
coincides with the strong-control expression
\eqref{eq:agreementGcontrol}. For general bounded Borel payoffs,
\eqref{eq:multT} is the canonical extension used throughout the paper.
\begin{proposition}
\label{prop:multDPP}
Let \(T_t\) be defined by \eqref{eq:multT}. Then, for every bounded
upper semianalytic function
\(\varphi:\mathbb R^d\to\mathbb R\), the following hold:
\begin{enumerate}[label=\textup{(\roman*)}]
 \item \(T_t\varphi\) is bounded and upper semianalytic for every
 \(t\ge0\);
 \item the semigroup identity
 \begin{equation}\label{eq:multsemigroup}
  T_{t+s}\varphi=T_t(T_s\varphi),
  \qquad s,t\ge0,
 \end{equation}
 holds pointwise on \(\mathbb R^d\);
 \item if \(t>0\) and \(\theta\le t\) is a bounded canonical stopping
 time, then
 \begin{equation}\label{eq:multstoppedDPP}
  T_t\varphi(x)
  =
  \sup_{P\in\mathcal P_{\mathcal Q}}
  E^P\!\left[
   T_{t-\theta}\varphi(X_\theta^x)
  \right],
  \qquad x\in\mathbb R^d;
 \end{equation}
 \item \(T_t\) is monotone, constant preserving, positively
 homogeneous, and subadditive.
\end{enumerate}
In \eqref{eq:multstoppedDPP}, the random variable inside the
expectation is interpreted with respect to the universal completion.
\end{proposition}
\begin{proof}
Assertion \textup{(i)} follows directly from
Lemma~\ref{lem:joint-value-measurability}. In particular,
\(T_s\varphi\) is again a bounded upper semianalytic function, so the
right-hand side of \eqref{eq:multsemigroup} is well defined.
For a finite canonical stopping time \(\tau\), define the admissible
continuation correspondence by
\[
 \Gamma_\tau(\omega):=\mathcal P_{\mathcal Q},
 \qquad \omega\in\Omega.
\]
Its graph is
\[
 \operatorname{Gr}(\Gamma_\tau)
 =
 \Omega\times\mathcal P_{\mathcal Q},
\]
up to the immaterial choice of the order of the two factors.
By Lemma~\ref{lem:joint-value-measurability},
\(\mathcal P_{\mathcal Q}\) is Borel; hence this graph is Borel and,
in particular, analytic.
Proposition~3.1 of Nutz--van Handel
\cite{NutzVanHandel2013} shows, for the constant compact convex
volatility set \(\mathcal Q\), that this family is invariant under
conditioning and stable under pasting. Thus the continuation family
satisfies Assumption~2.1 of that reference. Theorem~2.3 of
\cite{NutzVanHandel2013} therefore provides conditional upper
expectations, preserves upper semianalyticity, and gives the tower
property.
Fix \(s,t\ge0\), \(x\in\mathbb R^d\), and set
\[
 \xi(\omega):=\varphi(X_{t+s}^x(\omega)).
\]
By Lemma~\ref{lem:joint-value-measurability}, \(\xi\) is upper
semianalytic. Let \(\mathcal E_t\) denote the conditional upper
expectation at time \(t\). The conditional flow property in
Lemma~\ref{lem:universal-state-map} gives
\[
 X_{t+s}^x
   \bigl(\omega\otimes_t\widetilde\omega\bigr)
 =
 X_s^{X_t^x(\omega)}(\widetilde\omega)
\]
in the almost-sure sense required by conditioning and pasting.
Consequently,
\[
 \mathcal E_t(\xi)
 =
 T_s\varphi(X_t^x)
\]
under every \(P\in\mathcal P_{\mathcal Q}\). Applying the tower
property at times \(0\) and \(t\), we obtain
\[
\begin{aligned}
 T_{t+s}\varphi(x)
 &=\mathcal E_0(\xi)\\
 &=\mathcal E_0\bigl(\mathcal E_t(\xi)\bigr)\\
 &=\sup_{P\in\mathcal P_{\mathcal Q}}
   E^P\!\left[T_s\varphi(X_t^x)\right]\\
 &=T_t(T_s\varphi)(x).
\end{aligned}
\]
This proves \textup{(ii)}.
For \textup{(iii)}, fix \(t>0\) and a bounded canonical stopping time
\(\theta\le t\), and apply the same argument to
\[
 \xi(\omega):=\varphi(X_t^x(\omega)).
\]
The conditional flow property now gives
\[
 \mathcal E_\theta(\xi)
 =
 T_{t-\theta}\varphi(X_\theta^x)
\]
under every \(P\in\mathcal P_{\mathcal Q}\). By
Lemma~\ref{lem:joint-value-measurability} and
Remark~\ref{rem:stopped-value-measurability}, the right-hand side is
universally measurable. The stopping-time tower property in
Theorem~2.3 of \cite{NutzVanHandel2013} therefore yields
\[
\begin{aligned}
 T_t\varphi(x)
 &=\mathcal E_0(\xi)\\
 &=\mathcal E_0\bigl(\mathcal E_\theta(\xi)\bigr)\\
 &=\sup_{P\in\mathcal P_{\mathcal Q}}
   E^P\!\left[
    T_{t-\theta}\varphi(X_\theta^x)
   \right],
\end{aligned}
\]
which is \eqref{eq:multstoppedDPP}.
Finally, monotonicity, preservation of constants, positive
homogeneity, and subadditivity follow immediately from the
representation
\[
 T_t\varphi(x)
 =
 \sup_{P\in\mathcal P_{\mathcal Q}}
 E^P[\varphi(X_t^x)].
\]
\end{proof}
\begin{remark}[The chosen extension to Borel data]
The original $G$-expectation on $L_G^1(\Omega)$ does not by itself
determine values for all bounded Borel functionals. Equation
\eqref{eq:multT} specifies the extension used here through the
canonical covariance-law family $\mathcal P_{\Q}$. Its advantages are
that bounded upper semianalytic payoffs are admitted and that the DPP
follows from conditioning and pasting without selecting
state-dependent controls on a fixed Wiener space. The agreement
statement \eqref{eq:agreementGcontrol} is used only on the natural
$G$-expectation domain.
\end{remark}
\subsection{Bellman equations and positive-time regularisation}
\renewcommand{\theassumption}{G}
\begin{assumption}[Multiplicative $G$-SDE hypotheses]
\label{ass:multG}\
\begin{enumerate}[label=\textup{(G\arabic*)}]
 \item The maps $b$ and $\sigma$ are globally Lipschitz. Both have at
 most linear growth.
 \item There are $0<\underline q\le\overline q<\infty$ such that
 \begin{equation}\label{eq:Qell}
  \underline q I_m\le Q\le\overline q I_m,
  \qquad Q\in\Q.
 \end{equation}
 \item There are $c>0$ and $C_0\ge0$ such that
 \begin{equation}\label{eq:multGLyap}
  \sup_{Q\in\Q}
  \left\{
   2\langle x,b(x)\rangle
   +\tr\bigl(\sigma(x)Q\sigma(x)^{\mathsf T}\bigr)
  \right\}
  \le-c|x|^2+C_0.
 \end{equation}
 \item On every bounded set $O\subset\R^d$ there are
 $0<\lambda_O\le\Lambda_O<\infty$ such that
 \begin{equation}\label{eq:effectiveellmult}
  \lambda_OI_d
  \le\sigma(x)Q\sigma(x)^{\mathsf T}
  \le\Lambda_OI_d,
  \qquad x\in O,\quad Q\in\Q.
 \end{equation}
\end{enumerate}
\end{assumption}
If \eqref{eq:Qell} holds and $\sigma(x)$ has full row rank for every
$x$, then \eqref{eq:effectiveellmult} follows on compact sets.
For $Q\in\Q$, write
\[
 a_Q(x):=\sigma(x)Q\sigma(x)^{\mathsf T},
 \qquad
 \mathcal L^Qf(x)
 :=b(x)\cdot Df(x)+\frac12\tr\bigl(a_Q(x)D^2f(x)\bigr).
\]
Now we are in a position to derive upper and lower Bellman equations and establish the same Pucci class to 
verify the Harnack inequalities for viscosity solutions of both equations. 
\begin{proposition}
\label{prop:multbellman}
Suppose that Assumption~\ref{ass:multG} holds. Let
$\varphi:\R^d\to\R$ be bounded and uniformly continuous, and set
\[
 u(t,x):=T_t\varphi(x).
\]
Then $u$ is continuous on $[0,\infty)\times\R^d$ and is a viscosity
solution on $(0,\infty)\times\R^d$ of
\begin{equation}\label{eq:multBellman}
 \partial_tu
 =\sup_{Q\in\Q}
  \left\{
   b(x)\cdot Du
   +\frac12\tr\bigl(a_Q(x)D^2u\bigr)
  \right\}.
\end{equation}
If $M$ is constant and $w=M-u$, then $w$ is a continuous viscosity
solution of the dual equation
\begin{equation}\label{eq:multdual}
 \partial_tw
 =\inf_{Q\in\Q}
  \left\{
   b(x)\cdot Dw
   +\frac12\tr\bigl(a_Q(x)D^2w\bigr)
  \right\}.
\end{equation}
Let $O\subset\R^d$ be bounded and set
\[
 B_O:=\|b\|_{L^\infty(O)}.
\]
Every continuous viscosity solution $z$ of either equation belongs on
$(0,\infty)\times O$ to the same parabolic Pucci class:
\begin{equation}\label{eq:multPucciclass}
 \mathcal M^-_{\lambda_O/2,\Lambda_O/2}(D^2z)-B_O|Dz|
 \leq\partial_tz
 \leq
 \mathcal M^+_{\lambda_O/2,\Lambda_O/2}(D^2z)+B_O|Dz|.
\end{equation}
Consequently, if
\[
 \mathcal C_{2r}(t_0,x_0)
 :=(t_0-4r^2,t_0)\times B_{2r}(x_0)
 \Subset(0,\infty)\times O,
\]
then every nonnegative continuous viscosity solution of either
equation satisfies
\begin{equation}\label{eq:multlocalharnack}
 \sup_{\mathcal C_r^-}z
 \leq C_H\inf_{\mathcal C_r^+}z,
\end{equation}
where, for example,
\[
 \mathcal C_r^-
 :=(t_0-3r^2,t_0-2r^2)\times B_r(x_0),
 \qquad
 \mathcal C_r^+
 :=(t_0-r^2,t_0)\times B_r(x_0).
\]
The constant $C_H$ depends only on
\[
 d,\qquad \lambda_O,\qquad\Lambda_O,\qquad rB_O,
\]
and the fixed relative geometry of the cylinders. In particular, the
same constant applies to the upper and lower equations.
\end{proposition}
\begin{proof}
The standard stability estimates for the controlled state equation,
together with uniform continuity of $\varphi$, show that $u$ is
continuous. The DPP in Proposition~\ref{prop:multDPP} gives
\begin{equation}\label{eq:multDPP}
 u(t+s,x)
 =\sup_{P\in\mathcal P_{\Q}}
   E^P\bigl[u(t,X_s^x)\bigr],
 \qquad s,t\geq0.
\end{equation}
We recall the viscosity argument \cite{CIL1992,IshiiLions1990}. Fix
$(t_0,x_0)\in(0,\infty)\times\R^d$ and let
$\psi\in C^{1,2}$ touch $u$ from above at $(t_0,x_0)$. After the
standard localisation of the contact, choose an $h^2$-optimal law
$P_h\in\mathcal P_{\Q}$ in \eqref{eq:multDPP}. Then
\[
 u(t_0,x_0)
 \leq
 E^{P_h}\bigl[u(t_0-h,X_h^{x_0})\bigr]+h^2,
\]
and hence
\[
 0
 \leq
 E^{P_h}\left[
  \psi(t_0-h,X_h^{x_0})-\psi(t_0,x_0)
 \right]+o(h).
\]
Writing
\[
 d\langle B\rangle_s=\widehat Q_s^{P_h}\,ds,
 \qquad \widehat Q_s^{P_h}\in\Q,
\]
and applying It\^o's formula gives
\[
\begin{split}
0\leq
E^{P_h}\int_0^h
\biggl[
 &-\partial_t\psi(t_0-s,X_s^{x_0})
 +b(X_s^{x_0})\cdot D\psi(t_0-s,X_s^{x_0})\\
 &+\frac12\tr\left(
 a_{\widehat Q_s^{P_h}}(X_s^{x_0})
 D^2\psi(t_0-s,X_s^{x_0})
 \right)
\biggr]\,ds+o(h).
\end{split}
\]
Dividing by $h$ and letting $h\downarrow0$, using compactness of
$\Q$, yields
\[
 \partial_t\psi(t_0,x_0)
 \leq
 \sup_{Q\in\Q}
 \left\{
 b(x_0)\cdot D\psi(t_0,x_0)
 +\frac12\tr\bigl(a_Q(x_0)D^2\psi(t_0,x_0)\bigr)
 \right\}.
\]
Thus $u$ is a viscosity subsolution.
Suppose next that $\psi$ touches $u$ from below at $(t_0,x_0)$.
Fix an arbitrary $Q\in\Q$ and use in the DPP the law under which
$d\langle B\rangle_s=Q\,ds$. Then
\[
 u(t_0,x_0)
 \geq
 E^Q\bigl[u(t_0-h,X_h^{x_0})\bigr].
\]
Since $u\geq\psi$ near the contact point,
\[
 0
 \leq
 \psi(t_0,x_0)
 -
 E^Q\bigl[\psi(t_0-h,X_h^{x_0})\bigr]+o(h).
\]
It\^o's formula and $h\downarrow0$ give
\[
 \partial_t\psi(t_0,x_0)
 \geq
 b(x_0)\cdot D\psi(t_0,x_0)
 +\frac12\tr\bigl(a_Q(x_0)D^2\psi(t_0,x_0)\bigr).
\]
Since $Q$ was arbitrary, taking the supremum over $Q$ proves the
viscosity supersolution inequality. Hence $u$ solves
\eqref{eq:multBellman}.
Since every $\mathcal L^Q$ is linear and annihilates constants,
\[
 -\sup_{Q\in\Q}\mathcal L^Qu
 =
 \inf_{Q\in\Q}\mathcal L^Q(M-u).
\]
The standard change of test functions therefore shows that $w=M-u$
solves \eqref{eq:multdual}.
It remains to identify the structural class. By
\eqref{eq:effectiveellmult}, for $x\in O$, $A\in\mathbb S^d$ and
$p\in\R^d$,
\[
 \mathcal M^-_{\lambda_O/2,\Lambda_O/2}(A)-B_O|p|
 \leq
 b(x)\cdot p+\frac12\tr(a_Q(x)A)
 \leq
 \mathcal M^+_{\lambda_O/2,\Lambda_O/2}(A)+B_O|p|.
\]
Taking either the supremum or the infimum over $Q$ proves
\eqref{eq:multPucciclass}. The parabolic Krylov--Safonov Harnack
inequality for this Pucci class then gives
\eqref{eq:multlocalharnack}; see
\cite{Krylov1987,Wang1992}. The dependence of $C_H$ follows by
parabolic rescaling.
\end{proof}
The following positive-time regularisation result provides a key to extend viscosity solutions to the case with Borel initial data and ultimately obtain the existence and ergodicity of invariant expectations.  
\begin{corollary}
\label{cor:multregularisation}
Suppose that Assumption~\ref{ass:multG} holds. For every
$\varphi\in B_b(\R^d)$, the value function
\[
 u(t,x):=T_t\varphi(x)
\]
is locally H\"older continuous on
$(0,\infty)\times\R^d$.
More precisely, for every $0<\tau<T<\infty$ and every compact
$K\subset\R^d$, there are $\alpha\in(0,1)$ and
$C_{\tau,T,K}<\infty$ such that
\begin{equation}\label{eq:multHolderregularisation}
 |T_t\varphi(x)-T_s\varphi(y)|
 \leq C_{\tau,T,K}\|\varphi\|_\infty
 \bigl(|x-y|+|t-s|^{1/2}\bigr)^\alpha
\end{equation}
for all $s,t\in[\tau,T]$ and $x,y\in K$. Consequently,
\[
 T_t\bigl(B_b(\R^d)\bigr)\subset C(\R^d),
 \qquad t>0,
\]
and $(T_t)_{t\geq0}$ is a sublinear Markov semigroup on
$B_b(\R^d)$.
Moreover, $u$ is a continuous viscosity solution of
\eqref{eq:multBellman} on $(0,\infty)\times\R^d$. For every constant
$M$, the function $M-u$ is a continuous viscosity solution of
\eqref{eq:multdual}.
\end{corollary}
\begin{proof}
Fix
\[
 \overline{\mathcal C_{4r}(t_0,x_0)}
 \subset (0,\infty)\times O,
\]
where $O\subset\R^d$ is bounded. By
\eqref{eq:multstoppedDPP}, the bounded universally measurable
function
\[
 u(t,x):=T_t\varphi(x)
\]
satisfies the stopped dynamic programming identity required in
Lemma~\ref{lem:BorelDataSmoothing}.
For $P\in\mathcal P_{\Q}$, let
$\widehat Q_s^P\in\Q$ denote the quadratic-variation density of the
canonical process. Up to the exit time from $O$, the state process has
drift and covariance characteristics
\[
 \beta_s^P=b(X_s^x),
 \qquad
 A_s^P
 =
 \sigma(X_s^x)\widehat Q_s^P
 \sigma(X_s^x)^{\mathsf T}.
\]
Assumption~\ref{ass:multG} gives, uniformly in
$P\in\mathcal P_{\Q}$,
\[
 |\beta_s^P|
 \leq B_O:=\|b\|_{L^\infty(O)},
 \qquad
 \lambda_OI_d\leq A_s^P\leq\Lambda_OI_d.
\]
All hypotheses of Lemma~\ref{lem:BorelDataSmoothing} are therefore
satisfied. Since
\[
 |u(t,x)|\leq\|\varphi\|_\infty,
\]
that lemma yields
\[
 [u]_{C^{\alpha/2,\alpha}
       (\mathcal C_r(t_0,x_0))}
 \leq
 C r^{-\alpha}\|\varphi\|_\infty,
\]
where $C$ depends only on the local structural constants.
A finite covering of $[\tau,T]\times K$ gives
\eqref{eq:multHolderregularisation}. In particular,
$T_t\varphi$ is a continuous, hence Borel, function for every
$t>0$. Thus the semigroup identity of
Proposition~\ref{prop:multDPP} is an identity in $B_b(\R^d)$, proving
that $(T_t)_{t\geq0}$ is a sublinear Markov semigroup on
$B_b(\R^d)$.
It remains to identify the equation satisfied by $u$. Now that
continuity has been established, the test-function argument in the
proof of Proposition~\ref{prop:multbellman} applies directly. Indeed,
at an upper contact one chooses an $o(h)$-optimal law in the DPP. At a
lower contact one uses, successively, every law with constant
covariance density $Q\in\Q$. It\^o's formula gives respectively the
viscosity subsolution and supersolution inequalities for
\eqref{eq:multBellman}. This argument uses only the DPP, continuity of
$u$, and the local coefficient bounds; it does not require any
regularity of the initial function $\varphi$.
Thus $u$ is a continuous viscosity solution of
\eqref{eq:multBellman}. Since every $\mathcal L^Q$ is linear and
annihilates constants, the transformation $w=M-u$ gives a continuous
viscosity solution of \eqref{eq:multdual}.
\end{proof}
\subsection{Invariant expectation and exponential convergence}
\begin{theorem}
\label{thm:multG}
Suppose that Assumption~\ref{ass:multG} holds and let
$V(x)=|x|^2$. Then the operators \eqref{eq:multT} extend canonically
to a sublinear Markov semigroup on $B_V(\R^d)$. There exists a unique
invariant sublinear expectation $\Et\in\mathfrak E_V$. Moreover,
there are constants $C,\rho>0$ such that
\begin{equation}\label{eq:multGmix}
 |T_t\varphi(x)-\Et[\varphi]|
 \le Ce^{-\rho t}(1+|x|^2)\osc(\varphi),
 \qquad t\ge0,\quad\varphi\in B_b(\R^d).
\end{equation}
Consequently,
\[
 T_t\varphi(x)\longrightarrow\Et[\varphi]
\]
locally uniformly in $x$ for every $\varphi\in B_b(\R^d)$.
\end{theorem}
\begin{proof}
It\^o's formula under each $P\in\mathcal P_{\Q}$ and
\eqref{eq:multGLyap} give
\begin{equation}\label{eq:multLyapunovestimate}
 T_tV(x)
 \le e^{-ct}V(x)+\frac{C_0}{c}(1-e^{-ct}),
 \qquad t\ge0.
\end{equation}
Lemma~\ref{lem:BbBVbridge} therefore extends the canonical semigroup,
the DPP, and the positive-time regularity from $B_b(\R^d)$ to
$B_V(\R^d)$.
Fix $a>0$ and put
\[
 \gamma:=e^{-ca},
 \qquad
 K:=\frac{C_0}{c}(1-e^{-ca}).
\]
Then $T_aV\le\gamma V+K$, and
\eqref{eq:multLyapunovestimate} also gives the required short-time
Lyapunov bound. Choose
\[
 R>\frac{2K}{1-\gamma}=\frac{2C_0}{c},
 \qquad
 C_R:=\{x\in\R^d:V(x)\le R\}.
\]
Let $O_R$ be a bounded connected open set with
$C_R\Subset O_R$.
Choose $s_*\in(0,a)$ and $x_*\in C_R$. Compactness of $C_R$ and
connectedness of $O_R$ provide a finite parabolic Harnack chain from
$(s_*,x_*)$ to $\{a\}\times C_R$. Let
\[
 \mathfrak Q_R\Subset(0,\infty)\times O_R
\]
be a bounded positive-time set containing the cylinders in this chain
and $\{a\}\times C_R$.
For $\varphi\in B_V(\R^d)$, put $u(t,x)=T_t\varphi(x)$ and set
\[
 m:=\inf_{\mathfrak Q_R}u,
 \qquad
 M:=\sup_{\mathfrak Q_R}u.
\]
By Lemma~\ref{lem:BbBVbridge}, $u$ is bounded and continuous on
$\mathfrak Q_R$. Moreover, $u-m$ is a nonnegative solution of the
upper Bellman equation, while $M-u$ is a nonnegative solution of the
dual equation. Iterating \eqref{eq:multlocalharnack} along the finite
chain verifies the two-sided Harnack hypothesis of
Lemma~\ref{lem:Hosc}. Hence there is $\vartheta_R\in(0,1)$,
independent of $\varphi$, such that
\begin{equation}\label{eq:multlocaloscillation}
 \osc_{C_R}(T_a\varphi)
 \le\vartheta_R\osc_{\mathfrak Q_R}(T_\cdot\varphi),
 \qquad\varphi\in B_V(\R^d).
\end{equation}
Since $\mathfrak Q_R$ is a bounded positive-time set and its spatial
projection is contained in the bounded set $O_R$, there exist
$T_R\geq a$ and $\widehat R\geq R$ such that
\[
 \{a\}\times C_R
 \subset
 \mathfrak Q_R
 \subset
 (0,T_R]\times C_{\widehat R}.
\]
Thus Assumption~\ref{ass:LOC} is verified, while
Assumption~\ref{ass:L} follows from
\eqref{eq:multLyapunovestimate}. Therefore
Theorem~\ref{thm:weighted} and
Corollary~\ref{thm:invariant} yield the unique invariant
sublinear expectation and the exponential estimate
\eqref{eq:multGmix}.
\end{proof}
\section{The full \texorpdfstring{$G$}{G}-SDE}
\label{sec:fullG}
We now consider
\begin{equation}\label{eq:fullGSDE}
 dX_t=b(X_t)\,dt
 +\sum_{i,j=1}^m h_{ij}(X_t)\,d\langle B^i,B^j\rangle_t
 +\sum_{j=1}^m\sigma_j(X_t)\,dB_t^j.
\end{equation}
Here $X_t\in\R^d$, $h_{ij}(x)\in\R^d$, and the columns of
$\sigma(x)\in\R^{d\times m}$ are denoted by $\sigma_j(x)$.
The matrix-valued quadratic variation of the $G$-Brownian motion is
denoted by
\[
 \langle B\rangle_t
 :=
 \bigl(\langle B^i,B^j\rangle_t\bigr)_{i,j=1}^m.
\]
Its mutual quadratic-covariation entries may be defined by polarisation:
\[
 \langle B^i,B^j\rangle_t
 =
 \frac14\left(
  \langle B^i+B^j\rangle_t
  -\langle B^i-B^j\rangle_t
 \right).
\]
Unlike the quadratic variation of a Brownian motion with fixed covariance,
$\langle B\rangle$ is generally nondeterministic. Under every representing
law $P\in\mathcal P_{\Q}$, it agrees $P$-almost surely with the usual
quadratic variation of the canonical process and satisfies
\begin{equation}\label{eq:fullQVdensity}
 d\langle B\rangle_t=\widehat Q_t^P\,dt,
 \qquad
 \widehat Q_t^P\in\Q,
 \quad dt\otimes P\text{-a.e.}
\end{equation}
Consequently,
\[
 \sum_{i,j=1}^m
 h_{ij}(X_t)\,d\langle B^i,B^j\rangle_t
 =
 \sum_{i,j=1}^m
 (\widehat Q_t^P)^{ij}h_{ij}(X_t)\,dt,
 \qquad P\text{-a.s.}
\]
Thus the quadratic-variation term in \eqref{eq:fullGSDE} acts, under each
representing law, as an additional drift depending on the current
covariance scenario.
Since $\langle B\rangle$ is symmetric, only
\begin{equation}\label{eq:hsym}
 \widehat h_{ij}:=\frac12(h_{ij}+h_{ji})
\end{equation}
contributes to \eqref{eq:fullGSDE}. We may therefore assume that
$h_{ij}=h_{ji}$, although the formulas below are written so that no such
normalisation is required.
For $Q\in\Q$, define the effective drift and covariance by
\begin{align}
 b_Q(x)&=b(x)+\sum_{i,j=1}^mQ^{ij}h_{ij}(x),
 \label{eq:bQ}\\
 a_Q(x)&=\sigma(x)Q\sigma(x)^{\mathsf T}.
 \label{eq:aQ}
\end{align}
\subsection{Structural assumptions and the robust Lyapunov condition}
\renewcommand{\theassumption}{F}
\begin{assumption}[Full $G$-SDE hypotheses]\label{ass:fullG}\ 
\begin{enumerate}[label=\textup{(F\arabic*)}]
 \item $\Q$ is nonempty, compact and convex, and satisfies
 \eqref{eq:Qell}.
 \item $b,h_{ij},\sigma_j$ are globally Lipschitz and have at most linear
 growth.
 \item With $V(x)=\abs{x}^2$, there are $c>0$ and $C_0\geq0$ such that
 \begin{equation}\label{eq:fullLyap}
 \sup_{Q\in\Q}\left\{
  2\langle x,b(x)\rangle
  +2\sum_{i,j=1}^mQ^{ij}\langle x,h_{ij}(x)\rangle
  +\tr\bigl(\sigma(x)Q\sigma(x)^{\mathsf T}\bigr)
 \right\}
 \leq-c\abs{x}^2+C_0.
 \end{equation}
 \item On every bounded set $O\subset\R^d$ there are
 $0<\lambda_O\leq\Lambda_O<\infty$ such that
 \begin{equation}\label{eq:fullell}
  \lambda_OI_d\leq a_Q(x)\leq\Lambda_OI_d,
  \qquad x\in O,\quad Q\in\Q.
 \end{equation}
\end{enumerate}
\end{assumption}
No smoothness beyond (F2) is needed for the viscosity formulation used
below. Additional local smoothness may be imposed if classical solutions are
desired. The following provides the equivalent $G$-form of the Lyapunov condition.
\begin{proposition}
\label{prop:Glyapform}
Condition \eqref{eq:fullLyap} is equivalent to
\begin{equation}\label{eq:fullLyapG}
 2\langle x,b(x)\rangle+G(\mathcal H_V(x))
 \leq-c\abs{x}^2+C_0,
\end{equation}
where
\begin{equation}\label{eq:HV}
 (\mathcal H_V)_{ij}(x)
 =2\langle\sigma_i(x),\sigma_j(x)\rangle
 +2\langle x,h_{ij}(x)+h_{ji}(x)\rangle.
\end{equation}
For symmetric $h$, the second term is $4\langle x,h_{ij}(x)\rangle$.
\end{proposition}
\begin{proof}
By the definition of $G$,
\[
 G(\mathcal H_V(x))
 =
 \frac12\sup_{Q\in\Q}
 \sum_{i,j=1}^m Q^{ij}(\mathcal H_V)_{ij}(x).
\]
Using \eqref{eq:HV} and the symmetry of $Q$, we obtain
\[
 G(\mathcal H_V(x))
 =
 \sup_{Q\in\Q}\left\{
  \tr\bigl(\sigma(x)Q\sigma(x)^{\mathsf T}\bigr)
  +2\sum_{i,j=1}^mQ^{ij}\langle x,h_{ij}(x)\rangle
 \right\}.
\]
Thus \eqref{eq:fullLyapG} is exactly \eqref{eq:fullLyap}.
\end{proof}
\begin{lemma}
\label{lem:boundedhsigma}
Suppose that
\[
 \langle x,b(x)\rangle\leq-\lambda_b\abs{x}^2+C_b
\]
for some $\lambda_b>0$, and suppose that every $h_{ij}$ and $\sigma$ is
bounded. Then \eqref{eq:fullLyap} holds for suitable $c>0$ and
$C_0<\infty$.
\end{lemma}
\begin{proof}
Compactness of $\Q$ and boundedness of $h$ yield a constant $C_h$ such that
\[
 2\left|\sum_{i,j}Q^{ij}\langle x,h_{ij}(x)\rangle\right|
 \leq C_h\abs{x}
\]
uniformly in $Q$. Young's inequality absorbs this term into
$\varepsilon\abs{x}^2+C_\varepsilon$, while the covariance trace is
uniformly bounded. Taking $\varepsilon<2\lambda_b$ proves the claim.
\end{proof}
\begin{remark}
If $h$ or $\sigma$ has linear growth, its contribution to
\eqref{eq:fullLyap} may be quadratic. Dissipativity of $b$ alone is then
insufficient; the quadratic terms must satisfy a genuine balance condition.
Condition \eqref{eq:fullLyap} records this balance without discarding the
dependence on the covariance scenario.
\end{remark}
\subsection{Canonical representation and generator}
Let
\[
 \Omega:=C_0(\R_+;\R^m)
\]
be the canonical path space, let $B$ be the coordinate process, and let
$\mathbb F=(\mathcal F_t)_{t\geq0}$ be its raw filtration. Define
\[
 \mathcal P_{\Q}
 :=
 \left\{
 P\in\mathfrak P(\Omega):
 \begin{array}{l}
 B\text{ is a continuous local }P\text{-martingale,}\\[1mm]
 d\langle B\rangle_t=\widehat Q_t^P\,dt,
 \quad \widehat Q_t^P\in\Q,
 \quad dt\otimes P\text{-a.e.}
 \end{array}
 \right\}.
\]
The quadratic variation and its density are understood in their canonical
pathwise versions \cite{Nutz2012}. Condition~\textup{(F1)} of Assumption~\ref{ass:fullG} implies that this family is
nonempty. More importantly, Proposition~3.1 of
Nutz--van Handel~\cite{NutzVanHandel2013} verifies, for the constant compact
convex volatility set $\Q$, all three parts of their Assumption~2.1: the
graph of the conditional family is analytic, the family is invariant under
conditioning, and it is stable under pasting. Their Theorem~2.3 therefore
provides conditional upper expectations for upper semianalytic payoffs and
the tower property. These facts will be used explicitly in
Proposition~\ref{prop:fullDPP}.
Under $P\in\mathcal P_{\Q}$, uniform positive definiteness of $\Q$ allows
one to write
\begin{equation}\label{eq:fullcontrol}
 dX_t^x
 =b_{\widehat Q_t^P}(X_t^x)\,dt
 +\sigma(X_t^x)(\widehat Q_t^P)^{1/2}\,dW_t^P,
 \qquad X_0^x=x,
\end{equation}
for a $P$-Brownian motion $W^P$. Under condition~\textup{(F2)} of Assumption~\ref{ass:fullG},
the canonical $G$-SDE has a unique solution.
\begin{lemma}
\label{lem:full-universal-state}
Suppose that \(\mathcal Q\subset\mathbb S_+^m\) is compact and that
\(b,h_{ij}\) and \(\sigma\) are globally Lipschitz and have at most
linear growth. There exists a map
\[
        (x,\omega,t)\longmapsto X_t^x(\omega)
\]
which is
\[
\mathcal B(\mathbb R^d)\otimes\operatorname{Prog}(\mathbb F)
\text{-measurable},
\]
such that, for every \(P\in\mathcal P_{\mathcal Q}\) and every
\(x\in\mathbb R^d\),
\begin{equation}
\label{eq:full-universal-sde}
X_t^x
 =x+\int_0^t b(X_s^x)\,ds
   +\sum_{i,j=1}^m\int_0^t h_{ij}(X_s^x)\,
                   d\langle B^i,B^j\rangle_s
   +\int_0^t\sigma(X_s^x)\,dB_s ,
\qquad t\geq0,
\end{equation}
\(P\)-almost surely. In particular, \(X^x\) has continuous paths
\(P\)-almost surely and is indistinguishable from the unique
\(P\)-solution of \eqref{eq:full-universal-sde}.
Moreover, for every \(P\in\mathcal P_{\mathcal Q}\), every finite
stopping time \(\tau\), and every admissible continuation kernel
\(\nu(\omega)\in\mathcal P_{\mathcal Q}\), one has
\begin{equation}
\label{eq:full-conditional-flow}
 X_{\tau(\omega)+s}^x
       \bigl(\omega\otimes_\tau\widetilde\omega\bigr)
 =
 X_s^{\,X_\tau^x(\omega)}(\widetilde\omega),
 \qquad s\geq0,
\end{equation}
for
\[
        P(d\omega)\nu(\omega,d\widetilde\omega)
        \text{-almost every }(\omega,\widetilde\omega).
\]
Consequently, if \(\varphi:\mathbb R^d\to\mathbb R\) is bounded and
upper semianalytic, then
\[
 (t,x,P)\longmapsto E^P[\varphi(X_t^x)]
\]
is upper semianalytic. Hence
\[
 (t,x)\longmapsto
 T_t\varphi(x):=\sup_{P\in\mathcal P_{\mathcal Q}}
                       E^P[\varphi(X_t^x)]
\]
is jointly upper semianalytic.
\end{lemma}
\begin{proof}
We indicate only the difference from Lemma~\ref{lem:universal-state-map}.
Let
\[
       \Delta_n:=2^{-n},
       \qquad
       \eta_n(s):=\Delta_n
                   \left\lfloor\frac{s}{\Delta_n}\right\rfloor ,
\]
and define the Euler approximation by
\begin{align}
X_t^{n,x}
={}&x+\int_0^t b\bigl(X_{\eta_n(s)}^{n,x}\bigr)\,ds \notag\\
 &+\sum_{i,j=1}^m
      \int_0^t h_{ij}\bigl(X_{\eta_n(s)}^{n,x}\bigr)\,
                    d\langle B^i,B^j\rangle_s
   +\int_0^t\sigma\bigl(X_{\eta_n(s)}^{n,x}\bigr)\,dB_s .
\label{eq:full-euler}
\end{align}
Every integral in \eqref{eq:full-euler} is a finite sum of increments
of \(B\) and its pathwise quadratic variation on bounded time
intervals. Therefore
\[
        (x,\omega,t)\longmapsto X_t^{n,x}(\omega)
\]
is jointly Borel measurable and progressively measurable in
\((\omega,t)\).
Fix \(P\in\mathcal P_{\mathcal Q}\). Writing
\[
        d\langle B\rangle_s=\widehat Q_s^P\,ds,
        \qquad \widehat Q_s^P\in\mathcal Q,
        \quad ds\otimes P\text{-a.e.},
\]
equation \eqref{eq:full-universal-sde} becomes
\[
 dX_s^x
   =b_{\widehat Q_s^P}(X_s^x)\,ds+\sigma(X_s^x)\,dB_s.
\]
Compactness of \(\mathcal Q\) and the hypotheses on \(b\) and
\(h_{ij}\) imply that the maps \(b_Q\), \(Q\in\mathcal Q\), are
globally Lipschitz and have at most linear growth with constants
independent of \(Q\). The BDG and Gronwall estimates used in the proof
of Lemma~\ref{lem:universal-state-map} therefore give, for \(p\geq2\),
\(T>0\), and compact \(K\subset\mathbb R^d\),
\[
 \sup_{\substack{P\in\mathcal P_{\mathcal Q}\\x\in K}}
 E^P\!\left[
       \sup_{0\leq s\leq T}
       |X_s^{n,x}-X_s^{x,P}|^p
      \right]
 \leq C_{p,T,K}\Delta_n^{p/2},
\]
where \(X^{x,P}\) denotes the unique \(P\)-solution. Thus, for every
fixed \(P\) and \(x\), the Borel--Cantelli argument gives
\[
       X^{n,x}\longrightarrow X^{x,P}
       \quad\text{locally uniformly in time, }P\text{-a.s.}
\]
As in Lemma~\ref{lem:universal-state-map}, define \(X_t^x\) by the
coordinatewise pointwise \(\limsup\) of \(X_t^{n,x}\), setting it
equal to \(x\) if one of these limits is not finite. This definition
is
\(\mathcal B(\mathbb R^d)\otimes\Prog(\mathbb F)\)-measurable and
coincides with \(X^{x,P}\), \(P\)-almost surely, for every
\(P\in\mathcal P_{\mathcal Q}\) and every \(x\).
It remains only to note that conditioning and pasting preserve both
the canonical martingale part and the identity
\[
       d\langle B\rangle_s=\widehat Q_s\,ds,
       \qquad \widehat Q_s\in\mathcal Q.
\]
After concatenation at \(\tau\), both sides of
\eqref{eq:full-conditional-flow} therefore solve, under
\(\nu(\omega)\), the same equation with initial value
\(X_\tau^x(\omega)\), the same canonical noise, and the same
quadratic-variation term. Pathwise uniqueness yields
\eqref{eq:full-conditional-flow}.
Finally, the asserted upper semianalyticity follows from the joint
Borel measurability of the state map, the measurable-integration
theorem, and the analytic-graph property of
\(\mathcal P_{\mathcal Q}\).
\end{proof}
For every bounded upper semianalytic function $\varphi$, define
\begin{equation}\label{eq:fullT}
 T_t\varphi(x)
 :=
 \sup_{P\in\mathcal P_{\Q}}E^P\!\left[\varphi(X_t^x)\right].
\end{equation}
\begin{proposition}
\label{prop:fullgenerator}
For smooth test functions, the generator of \eqref{eq:fullT} is
\begin{equation}\label{eq:Fcontrol}
 F(x,p,A)=\sup_{Q\in\Q}\left\{
 b_Q(x)\cdot p+\frac12\tr(a_Q(x)A)
 \right\}.
\end{equation}
Equivalently,
\begin{equation}\label{eq:FG}
 F(x,p,A)=b(x)\cdot p+G(H(x,p,A)),
\end{equation}
where $H(x,p,A)\in\mathbb S^m$ is given by
\begin{equation}\label{eq:Hmatrix}
 H_{ij}(x,p,A)
 =\langle A\sigma_i(x),\sigma_j(x)\rangle
 +\langle p,h_{ij}(x)+h_{ji}(x)\rangle.
\end{equation}
If $h_{ij}=h_{ji}$, the last term is $2\langle p,h_{ij}(x)\rangle$.
\end{proposition}
\begin{proof}
For a fixed $Q\in\Q$, It\^o's formula gives
\[
 \cL^Qf(x)
 =b_Q(x)\cdot Df(x)
 +\frac12\tr\bigl(a_Q(x)D^2f(x)\bigr).
\]
Taking the supremum yields \eqref{eq:Fcontrol}. Since $Q$ is symmetric,
\[
 \sum_{i,j}Q^{ij}\langle p,h_{ij}\rangle
 =\frac12\sum_{i,j}Q^{ij}\langle p,h_{ij}+h_{ji}\rangle.
\]
The definition of $G$ then gives \eqref{eq:FG}--\eqref{eq:Hmatrix}.
\end{proof}
\begin{remark}
\label{rem:factors}
With the convention $G(A)=\frac12\sup_{Q\in\Q}\tr(QA)$, the argument of
$G$ must contain the full Hessian term and the symmetrised $h$-term appearing
in \eqref{eq:Hmatrix}. Formula \eqref{eq:Fcontrol} is independent of this
notational convention and is therefore the safest form to use.
\end{remark}
\begin{remark}
\label{rem:fullcontrolBorel}
Under $P\in\mathcal P_{\Q}$, equation \eqref{eq:fullGSDE} reads
\begin{equation}\label{eq:fullcontrolrestated}
 X_t^x
 =x+\int_0^t b_{\widehat Q_s^P}(X_s^x)\,ds
   +\int_0^t\sigma(X_s^x)\,dB_s,
 \qquad P\text{-a.s.}
\end{equation}
Thus the uncertain covariance enters both the drift and the second-order
coefficient. This is the principal structural difference from the
multiplicative $G$-SDE considered in Section~\ref{sec:multG}.
By the canonical representation of $G$-expectation
\cite[Section~3]{NutzVanHandel2013},
\[
 \sup_{P\in\mathcal P_{\Q}}E^P[\xi]
 =
 \widehat{\mathbb E}[\xi],
 \qquad
 \xi\in L_G^1(\Omega_t).
\]
Consequently, whenever
$\varphi(X_t^x)\in L_G^1(\Omega_t)$, equation
\eqref{eq:fullT} satisfies
\[
 T_t\varphi(x)
 =
 \widehat{\mathbb E}\bigl[\varphi(X_t^x)\bigr].
\]
Thus the canonical-law and $G$-expectation semigroups agree on their
common domain, in particular for
$\varphi\in C_{b,\mathrm{Lip}}(\R^d)$. For general bounded Borel
$\varphi$, equation \eqref{eq:fullT} is the canonical extension used
in this paper. This distinction is relevant because the original
functional on $L_G^1$ does not, by itself, determine a unique extension
to all Borel payoffs.
\end{remark}
\subsection{Dynamic programming and positive-time regularity}
\begin{proposition}
\label{prop:fullDPP}
Suppose that Assumption~\ref{ass:fullG} holds. Then, for every bounded
upper semianalytic function $\varphi:\R^d\to\R$, the following hold:
\begin{enumerate}[label=\textup{(\roman*)}]
 \item $T_t\varphi$ is bounded and upper semianalytic for every
 $t\geq0$;
 \item the semigroup identity
 \begin{equation}\label{eq:fullsemigroup}
  T_{t+s}\varphi=T_t(T_s\varphi),
  \qquad s,t\geq0,
 \end{equation}
 holds pointwise on $\R^d$;
 \item if $t>0$ and $\theta\leq t$ is a bounded canonical stopping
 time, then
 \begin{equation}\label{eq:fullstoppedDPP}
  T_t\varphi(x)
  =
  \sup_{P\in\mathcal P_{\Q}}
  E^P\!\left[
    T_{t-\theta}\varphi(X_\theta^x)
  \right],
  \qquad x\in\R^d;
 \end{equation}
 \item $T_t$ is monotone, constant preserving, positively homogeneous,
 and subadditive.
\end{enumerate}
In \eqref{eq:fullstoppedDPP}, the random variable inside the expectation
is interpreted with respect to the universal completion of
$\mathcal F_\theta$.
\end{proposition}
\begin{proof}
By Proposition~3.1 of Nutz--van Handel~\cite{NutzVanHandel2013},
the family \(\mathcal P_{\mathcal Q}\) has an analytic graph, is
invariant under conditioning, and is stable under pasting. Thus it
satisfies Assumption~2.1 of that reference. Theorem~2.3 therein
therefore provides conditional upper expectations that preserve
upper semianalyticity and satisfy the tower property.
Let \(\varphi:\mathbb R^d\to\mathbb R\) be bounded and upper
semianalytic. By Lemma~\ref{lem:full-universal-state}, the functional
\[
       (x,\omega)\longmapsto\varphi(X_{t+s}^x(\omega))
\]
is upper semianalytic. Moreover, the conditional flow identity
\eqref{eq:full-conditional-flow} implies that its conditional upper
expectation at time \(t\) is
\[
       T_s\varphi(X_t^x).
\]
Applying the tower property at times \(0\) and \(t\) gives
\[
\begin{aligned}
T_{t+s}\varphi(x)
  &=\sup_{P\in\mathcal P_{\mathcal Q}}
        E^P[\varphi(X_{t+s}^x)]  \\
  &=\sup_{P\in\mathcal P_{\mathcal Q}}
        E^P[T_s\varphi(X_t^x)]
   =T_t(T_s\varphi)(x).
\end{aligned}
\]
The joint upper semianalyticity conclusion of
Lemma~\ref{lem:full-universal-state} shows, in particular, that
\(T_s\varphi\) is upper semianalytic, so every term above is
well-defined on the universal completion.
The same argument at a bounded canonical stopping time
\(\theta\leq t\) yields
\[
       T_t\varphi(x)
       =
       \sup_{P\in\mathcal P_{\mathcal Q}}
       E^P\!\left[
          T_{t-\theta}\varphi(X_\theta^x)
       \right].
\]
Here the joint upper semianalyticity of
\[
       (s,y)\longmapsto T_s\varphi(y)
\]
ensures that
\(T_{t-\theta}\varphi(X_\theta^x)\) is universally measurable.
Finally, monotonicity, preservation of constants, positive
homogeneity, and subadditivity follow directly from the
representation
\[
       T_t\varphi(x)
       =
       \sup_{P\in\mathcal P_{\mathcal Q}}
       E^P[\varphi(X_t^x)].
\]
\end{proof}
\begin{proposition}
\label{prop:fullHarnack}
Suppose that Assumption~\ref{ass:fullG} holds. Let
$\varphi:\R^d\to\R$ be bounded and uniformly continuous, and set
\[
 u(t,x):=T_t\varphi(x).
\]
Then $u$ is continuous on $[0,\infty)\times\R^d$ and is a viscosity
solution on $(0,\infty)\times\R^d$ of
\begin{equation}\label{eq:fullupperPDE}
 \partial_tu
 =\sup_{Q\in\Q}
 \left\{
 b_Q(x)\cdot Du+\frac12\tr(a_Q(x)D^2u)
 \right\}.
\end{equation}
If $M$ is constant and $w=M-u$, then $w$ is a continuous viscosity
solution of
\begin{equation}\label{eq:fulllowerPDE}
 \partial_tw
 =\inf_{Q\in\Q}
 \left\{
 b_Q(x)\cdot Dw+\frac12\tr(a_Q(x)D^2w)
 \right\}.
\end{equation}
Fix a bounded open set $O\subset\R^d$ and define
\begin{equation}\label{eq:fullBO}
 B_O
 :=\sup_{\substack{x\in O\\Q\in\Q}}\abs{b_Q(x)}
 <\infty.
\end{equation}
Every continuous viscosity solution $z$ of either equation belongs on
$(0,\infty)\times O$ to the same parabolic Pucci class:
\begin{equation}\label{eq:fullPucciclass}
 \mathcal M^-_{\lambda_O/2,\Lambda_O/2}(D^2z)-B_O\abs{Dz}
 \leq \partial_tz
 \leq
 \mathcal M^+_{\lambda_O/2,\Lambda_O/2}(D^2z)+B_O\abs{Dz}.
\end{equation}
Consequently, nonnegative continuous viscosity solutions of either
equation satisfy the same local parabolic Harnack inequality.
\end{proposition}
\begin{proof}
The uniform stability estimates for the full controlled state
equation, together with uniform continuity of $\varphi$, show that
$u$ is continuous. The dynamic programming identity of
Proposition~\ref{prop:fullDPP} gives
\begin{equation}\label{eq:fullDPP}
 u(t+s,x)
 =\sup_{P\in\mathcal P_{\Q}}
 E^P\!\left[u(t,X_s^x)\right],
 \qquad s,t\geq0,
\end{equation}
and hence, for $0<h<t_0$,
\begin{equation}\label{eq:fullDPPbackward}
 u(t_0,x_0)
 =\sup_{P\in\mathcal P_{\Q}}
 E^P\!\left[u(t_0-h,X_h^{x_0})\right].
\end{equation}
Fix $(t_0,x_0)\in(0,\infty)\times\R^d$. Suppose that
$\psi\in C^{1,2}$ touches $u$ from above at $(t_0,x_0)$. After the
standard localisation of the contact, choose an $h^2$-optimal law
$P_h$ in \eqref{eq:fullDPPbackward}. Under $P_h$,
\[
 dX_s^{x_0}
 =
 b_{\widehat Q_s^{P_h}}(X_s^{x_0})\,ds
 +\sigma(X_s^{x_0})\,dB_s,
 \qquad
 d\langle B\rangle_s=\widehat Q_s^{P_h}\,ds.
\]
The DPP and It\^o's formula give
\[
\begin{split}
0\leq
E^{P_h}\int_0^h
\biggl[
 &-\partial_t\psi(t_0-s,X_s^{x_0})
 +b_{\widehat Q_s^{P_h}}(X_s^{x_0})
  \cdot D\psi(t_0-s,X_s^{x_0})\\
 &+\frac12\tr\left(
 a_{\widehat Q_s^{P_h}}(X_s^{x_0})
 D^2\psi(t_0-s,X_s^{x_0})
 \right)
\biggr]\,ds+o(h).
\end{split}
\]
Dividing by $h$ and letting $h\downarrow0$ yields
\[
 \partial_t\psi(t_0,x_0)
 \leq
 \sup_{Q\in\Q}
 \left\{
 b_Q(x_0)\cdot D\psi(t_0,x_0)
 +\frac12\tr\bigl(a_Q(x_0)D^2\psi(t_0,x_0)\bigr)
 \right\}.
\]
Thus $u$ is a viscosity subsolution.
Suppose next that $\psi$ touches $u$ from below at $(t_0,x_0)$.
Fix an arbitrary $Q\in\Q$ and use in
\eqref{eq:fullDPPbackward} the law having constant
quadratic-variation density $Q$. Since
\[
 u(t_0,x_0)
 \geq
 E^Q\bigl[u(t_0-h,X_h^{x_0})\bigr],
\]
the lower-contact inequality and It\^o's formula give
\[
 \partial_t\psi(t_0,x_0)
 \geq
 b_Q(x_0)\cdot D\psi(t_0,x_0)
 +\frac12\tr\bigl(a_Q(x_0)D^2\psi(t_0,x_0)\bigr).
\]
Since $Q$ is arbitrary, taking the supremum proves the viscosity
supersolution inequality. Hence $u$ solves
\eqref{eq:fullupperPDE}.
Since each $\cL^Q$ is linear and annihilates constants,
\[
 -\sup_{Q\in\Q}\cL^Qu
 =
 \inf_{Q\in\Q}\cL^Q(M-u).
\]
The standard change of test functions therefore proves
\eqref{eq:fulllowerPDE}.
It remains to identify the common structural class. Compactness of
$\Q$ and local boundedness of $b$ and $h$ imply \eqref{eq:fullBO};
more explicitly,
\[
 B_O
 \leq \|b\|_{L^\infty(O)}
 +\left(
   \sup_{Q\in\Q}\sum_{i,j=1}^m\abs{Q^{ij}}
  \right)
  \max_{1\leq i,j\leq m}\|h_{ij}\|_{L^\infty(O)}.
\]
Together with \eqref{eq:fullell}, this gives, for every
$A\in\mathbb S^d$ and $p\in\R^d$,
\begin{equation}\label{eq:fullPuccibound}
 \mathcal M^-_{\lambda_O/2,\Lambda_O/2}(A)-B_O\abs{p}
 \leq
 b_Q(x)\cdot p+\frac12\tr(a_Q(x)A)
 \leq
 \mathcal M^+_{\lambda_O/2,\Lambda_O/2}(A)+B_O\abs{p}.
\end{equation}
Taking either the supremum or the infimum over $Q$ proves
\eqref{eq:fullPucciclass}.
If
\[
 \mathcal C_{2r}(t_0,x_0)
 :=(t_0-4r^2,t_0)\times B_{2r}(x_0)
 \Subset(0,\infty)\times O,
\]
the parabolic Krylov--Safonov Harnack theorem applied to
\eqref{eq:fullPucciclass} yields
\begin{equation}\label{eq:fullHarnackexplicit}
 \sup_{\mathcal C_r^-}z
 \leq C_H\inf_{\mathcal C_r^+}z
\end{equation}
for every nonnegative continuous viscosity solution $z$ of either
equation. Here one may take
\[
 \mathcal C_r^-
 =(t_0-3r^2,t_0-2r^2)\times B_r(x_0),
 \qquad
 \mathcal C_r^+
 =(t_0-r^2,t_0)\times B_r(x_0).
\]
The constant depends only on
$d$, $\lambda_O$, $\Lambda_O$, $rB_O$ and the fixed relative geometry
of the cylinders. In particular, it is independent of
$z$, $M$, $\varphi$ and the covariance scenario.
\end{proof}
\begin{corollary}
\label{cor:fullregularisation}
Suppose that Assumption~\ref{ass:fullG} holds. For every
$\varphi\in B_b(\R^d)$, the value function
\[
 u(t,x):=T_t\varphi(x)
\]
is locally H\"older continuous on $(0,\infty)\times\R^d$.
More precisely, for every $0<\tau<T<\infty$ and every compact
$K\subset\R^d$, there are $\alpha\in(0,1)$ and
$C_{\tau,T,K}<\infty$ such that
\begin{equation}\label{eq:fullHolderregularisation}
 \abs{T_t\varphi(x)-T_s\varphi(y)}
 \leq C_{\tau,T,K}\|\varphi\|_\infty
 \left(\abs{x-y}+\abs{t-s}^{1/2}\right)^\alpha
\end{equation}
for all $s,t\in[\tau,T]$ and $x,y\in K$. Consequently,
\[
 T_t\bigl(B_b(\R^d)\bigr)\subset C(\R^d),
 \qquad t>0,
\]
and $(T_t)_{t\geq0}$ is a sublinear Markov semigroup on
$B_b(\R^d)$.
Moreover, $u$ is a continuous viscosity solution of
\eqref{eq:fullupperPDE} on $(0,\infty)\times\R^d$. For every constant
$M$, the function $M-u$ is a continuous viscosity solution of
\eqref{eq:fulllowerPDE}.
\end{corollary}
\begin{proof}
Fix
\[
 \overline{\mathcal C_{4r}(t_0,x_0)}
 \subset (0,\infty)\times O,
\]
where $O\subset\R^d$ is bounded. By the stopping-time dynamic
programming identity in Proposition~\ref{prop:fullDPP}, the bounded
universally measurable value function
\[
 u(t,x):=T_t\varphi(x)
\]
satisfies the stopped DPP required in
Lemma~\ref{lem:BorelDataSmoothing}.
Under $P\in\mathcal P_{\Q}$, its local drift and covariance
characteristics are
\[
 \beta_s^P
 =
 b_{\widehat Q_s^P}(X_s^x),
 \qquad
 A_s^P
 =
 \sigma(X_s^x)\widehat Q_s^P
 \sigma(X_s^x)^{\mathsf T}.
\]
By Assumption~\ref{ass:fullG},
\[
 |\beta_s^P|
 \leq
 B_O
 :=
 \sup_{\substack{x\in O\\Q\in\Q}}\abs{b_Q(x)}
 <\infty,
 \qquad
 \lambda_OI_d\leq A_s^P\leq\Lambda_OI_d,
\]
uniformly in $P$. Lemma~\ref{lem:BorelDataSmoothing} therefore gives
\[
 [u]_{C^{\alpha/2,\alpha}
       (\mathcal C_r(t_0,x_0))}
 \leq
 C r^{-\alpha}\|\varphi\|_\infty.
\]
A finite covering of $[\tau,T]\times K$ proves
\eqref{eq:fullHolderregularisation}. In particular,
$T_t\varphi\in C(\R^d)$ for every $t>0$.
The stopped and deterministic dynamic programming identities in
Proposition~\ref{prop:fullDPP} now hold between bounded Borel
functions. Hence $(T_t)_{t\geq0}$ is a sublinear Markov semigroup on
$B_b(\R^d)$.
It remains to identify the equation satisfied by $u$. Since $u$ is
now continuous, the test-function argument in the proof of
Proposition~\ref{prop:fullHarnack} applies directly. At an upper
contact one uses an $o(h)$-optimal law, whose quadratic-variation
density may vary progressively. At a lower contact one uses,
successively, every law with fixed constant covariance
$Q\in\Q$. The DPP and It\^o's formula give the viscosity subsolution
and supersolution inequalities for \eqref{eq:fullupperPDE}. Thus $u$
is a continuous viscosity solution of the upper Bellman equation,
without any appeal to semicontinuous envelopes or measurable-data
viscosity theory.
Since each $\cL^Q$ is linear and annihilates constants, the
transformation $w=M-u$ gives a continuous viscosity solution of
\eqref{eq:fulllowerPDE}.
\end{proof}
\subsection{Invariant expectation and exponential convergence}
\begin{theorem}
\label{thm:fullG}
Suppose that Assumption~\ref{ass:fullG} holds and let
$V(x)=\abs{x}^2$. Then the operators \eqref{eq:fullT} extend canonically to
a sublinear Markov semigroup on $B_V(\R^d)$. There exists a unique invariant
sublinear expectation $\Et\in\mathfrak E_V$. Moreover, there are constants
$C,\rho>0$ such that
\begin{equation}\label{eq:fullGmix}
 \abs{T_t\varphi(x)-\Et[\varphi]}
 \leq Ce^{-\rho t}(1+\abs{x}^2)\osc(\varphi),
 \qquad t\geq0,\quad\varphi\in B_b(\R^d).
\end{equation}
Consequently, $T_t\varphi(x)\to\Et[\varphi]$ locally uniformly in $x$ for
every $\varphi\in B_b(\R^d)$.
\end{theorem}
\begin{proof}
It\^o's formula under each $P\in\mathcal P_{\Q}$, followed by
\eqref{eq:fullLyap}, gives the uniform estimate
\begin{equation}\label{eq:fullVbound}
 T_tV(x)
 \leq e^{-ct}V(x)+\frac{C_0}{c}(1-e^{-ct}),
 \qquad t\geq0.
\end{equation}
Lemma~\ref{lem:BbBVbridge} therefore extends the semigroup and the
positive-time regularity from $B_b(\R^d)$ to $B_V(\R^d)$.
Fix $a>0$ and set
\[
 \gamma=e^{-ca},
 \qquad
 K=\frac{C_0}{c}(1-e^{-ca}).
\]
Then $T_aV\leq\gamma V+K$, and \eqref{eq:fullVbound} also supplies the
required short-time Lyapunov bound. Choose
\[
 R>\frac{2K}{1-\gamma}=\frac{2C_0}{c},
 \qquad
 C_R:=\{x\in\R^d:V(x)\leq R\},
\]
and a bounded connected open set $O_R$ with $C_R\Subset O_R$. By
condition~\textup{(F4)} of Assumption~\ref{ass:fullG} and \eqref{eq:fullBO}, the upper and lower
Bellman equations belong on $O_R$ to one uniformly parabolic Pucci class,
with constants independent of the covariance scenario.
Fix $s_*\in(0,a)$ and $x_*\in C_R$. Compactness of $C_R$ and connectedness
of $O_R$ yield a finite parabolic Harnack chain, contained in
$(0,\infty)\times O_R$, joining $(s_*,x_*)$ to every point of
$\{a\}\times C_R$. Let $\mathfrak Q_R\Subset(0,\infty)\times O_R$ contain
all cylinders in this finite chain. For $\varphi\in B_V(\R^d)$, put
$u(t,x)=T_t\varphi(x)$. By Lemma~\ref{lem:BbBVbridge} and
Corollary~\ref{cor:fullregularisation}, $u$ is bounded and continuous on
$\mathfrak Q_R$. If
\[
 m=\inf_{\mathfrak Q_R}u,
 \qquad
 M=\sup_{\mathfrak Q_R}u,
\]
then $u-m$ is a nonnegative solution of the upper Bellman equation, whereas
$M-u$ is a nonnegative solution of the lower Bellman equation. Iterating the
common Harnack inequality \eqref{eq:fullHarnackexplicit} along the chain and
using Lemma~\ref{lem:Hosc} gives a number $\vartheta_R\in(0,1)$, independent
of $\varphi$, such that
\[
 \operatorname{osc}_{C_R}(T_a\varphi)
 \leq\vartheta_R
 \operatorname{osc}_{\mathfrak Q_R}(T_\cdot\varphi).
\]
This gives the oscillation estimate required in
Assumption~\ref{ass:LOC}. Since $\mathfrak Q_R$ is a bounded
positive-time set and its spatial projection is contained in $O_R$,
there exist $T_R\geq a$ and $\widehat R\geq R$ such that
\[
 \{a\}\times C_R
 \subset
 \mathfrak Q_R
 \subset
 (0,T_R]\times C_{\widehat R}.
\]
Thus Assumption~\ref{ass:LOC} is verified, while
Assumption~\ref{ass:L} follows from \eqref{eq:fullVbound}.
Theorem~\ref{thm:weighted} and
Corollary~\ref{thm:invariant} therefore give the unique invariant
sublinear expectation and \eqref{eq:fullGmix}.
\end{proof}
\begin{remark}
Condition~\textup{(F4)} of Assumption~\ref{ass:fullG} on every bounded set is used in
Proposition~\ref{prop:fullHarnack} and
Corollary~\ref{cor:fullregularisation} to construct a positive-time regular
semigroup on $B_b(\R^d)$, and hence, by
Lemma~\ref{lem:BbBVbridge}, on $B_V(\R^d)$. Once that framework is available,
the contraction step in Theorem~\ref{thm:fullG} uses ellipticity only on the
single bounded connected neighbourhood $O_R$ of the Lyapunov sublevel set
$C_R$. Thus, if dynamic programming, the $B_V$-extension and the required
positive-time regularity are obtained by another argument, (F4) may be
replaced for the contraction theorem itself by uniform ellipticity on such
an $O_R$.
A global bound
\[
 \sigma(x)Q\sigma(x)^{\mathsf T}\geq\lambda I_d,
 \qquad x\in\R^d,\quad Q\in\Q,
\]
is therefore convenient but stronger than the Harnack-contraction step
requires. If $\sigma$ is continuous and has full row rank everywhere, the
uniform covariance bounds \eqref{eq:Qell} imply \eqref{eq:fullell} on every
compact set.
\end{remark}
\section{Ergodicity}
\label{sec:ergodicity}
Feng and Zhao~\cite{FengZhao2021} associate with an invariant sublinear
expectation a canonical stationary sublinear dynamical system on path space.
We verify here that their ergodicity criterion applies to the semigroups
constructed in Sections~\ref{sec:multG} and~\ref{sec:fullG}.
\begin{definition}
\label{def:ergodicity}
Let $(\widehat\Omega,\widehat{\mathcal F},\widehat{\mathbb E})$ carry an
expectation-preserving measurable flow
$(\widehat\theta_t)_{t\in\R}$. The system is \emph{ergodic} if every invariant
set $A\in\widehat{\mathcal F}$, that is,
$\widehat\theta_t^{-1}A=A$ for every $t\in\R$, satisfies
\begin{equation}\label{eq:ergodicdef}
 \widehat{\mathbb E}[I_A]=0
 \quad\text{or}\quad
 \widehat{\mathbb E}[I_{A^c}]=0.
\end{equation}
An invariant sublinear expectation for a Markov semigroup is called
ergodic when its associated canonical stationary system is ergodic.
\end{definition}
\begin{proposition}
\label{prop:fixed}
Under the hypotheses of Corollary~\ref{thm:invariant}, suppose that
$\varphi\in B_b(E)$ satisfies
\begin{equation}\label{eq:fixed}
 T_t\varphi=\varphi,\qquad t\geq0.
\end{equation}
Then $\varphi$ is constant.
\end{proposition}
\begin{proof}
For every $x\in E$, estimate \eqref{eq:continuousrate} and
\eqref{eq:fixed} give
\[
 \varphi(x)
 =\lim_{t\to\infty}T_t\varphi(x)
 =\Et[\varphi].
\]
The right-hand side is independent of $x$.
\end{proof}
\begin{corollary}
\label{cor:ergodic}
For each of the sublinear Markov semigroups considered in
Theorems~\ref{thm:multG} and \ref{thm:fullG}, let $\Et$ be its unique
invariant sublinear expectation. Then the canonical stationary sublinear
dynamical system associated with $(T_t,\Et)$ is ergodic in the sense of
Definition~\ref{def:ergodicity}.
\end{corollary}
\begin{proof}
Since
\[
        B_b(\mathbb R^d)\subset B_V(\mathbb R^d),
\]
the restriction of \(\widetilde{\mathbb E}\) to
\(B_b(\mathbb R^d)\) is an invariant sublinear expectation for the
corresponding sublinear Markov semigroup. The nonlinear Kolmogorov
construction of Feng and Zhao therefore gives a canonical stationary
sublinear dynamical system
\[
 \bigl(
   \widehat\Omega,\widehat{\mathcal D},
   (\widehat\theta_t)_{t\in\mathbb R},
   \widehat{\mathbb E}_{\widetilde{\mathbb E}}
 \bigr)
\]
with canonical process
\[
        \widehat X_t(\widehat\omega)=\widehat\omega(t),
        \qquad t\in\mathbb R.
\]
For notational simplicity, write
\[
        \widehat{\mathbb E}
        :=\widehat{\mathbb E}_{\widetilde{\mathbb E}}.
\]
We first prove that \(\widehat X\) is stochastically continuous.
In both cases under consideration, the controlled drift and diffusion
coefficients are globally Lipschitz and have at most linear growth,
uniformly over \(P\in\mathcal P_{\mathcal Q}\). In the full
\(G\)-SDE case, the relevant drift under a covariance scenario
\(Q\in\mathcal Q\) is
\[
        b_Q(x)
        =
        b(x)+\sum_{i,j=1}^m Q_{ij}h_{ij}(x),
\]
which has uniform Lipschitz and linear-growth bounds because
\(\mathcal Q\) is compact.
Standard BDG and Gronwall estimates therefore give a constant
\(C<\infty\) such that, for \(0<h\leq1\),
\[
 \sup_{P\in\mathcal P_{\mathcal Q}}
 E^P\left[
      \sup_{0\leq r\leq h}|X_r^x-x|^2
     \right]
 \leq C(1+|x|^2)h.
\]
The same estimates give
\[
 \sup_{P\in\mathcal P_{\mathcal Q}}
 E^P\left[
      \sup_{0\leq r\leq1}|X_r^x-X_r^y|^2
     \right]
 \leq C|x-y|^2.
\]
Fix \(\delta>0\), and define
\[
 g_\delta(z)
 =
 \min\left\{
       1,\,
       \frac{2}{\delta}
       \max\left\{0,|z|-\frac{\delta}{2}\right\}
      \right\},
 \qquad z\in\mathbb R^d.
\]
Then \(g_\delta\) is bounded and Lipschitz, and
\[
        I_{\{|z|\geq\delta\}}
        \leq g_\delta(z)
        \leq \frac{4}{\delta^2}|z|^2.
\]
For \(0<h\leq1\), set
\[
 \Phi_h^\delta(x)
 =
 \sup_{P\in\mathcal P_{\mathcal Q}}
 E^P\left[g_\delta(X_h^x-x)\right].
\]
The measurable state-map constructions in Sections~\ref{sec:multG} and~\ref{sec:fullG} ensure
that this function is well defined. Moreover, if \(L_\delta\) denotes
the Lipschitz constant of \(g_\delta\), then
\[
\begin{split}
 |\Phi_h^\delta(x)-\Phi_h^\delta(y)|
 &\leq
 L_\delta
 \sup_{P\in\mathcal P_{\mathcal Q}}
 E^P\left[
      |X_h^x-X_h^y|+|x-y|
     \right]  \\
 &\leq C_\delta |x-y|.
\end{split}
\]
Consequently,
\[
        \Phi_h^\delta\in C_b(\mathbb R^d).
\]
The preceding small-time moment estimate also gives
\[
 0\leq\Phi_h^\delta(x)
 \leq
 \frac{4C}{\delta^2}(1+|x|^2)h
 =
 \frac{4C}{\delta^2}(1+V(x))h.
\]
Let \(s<t\) and put \(h=t-s\leq1\). By monotonicity, the definition
of the canonical two-time distribution, and stationarity,
\[
\begin{split}
 \widehat{\mathbb E}
 \left[
  I_{\{|\widehat X_t-\widehat X_s|\geq\delta\}}
 \right]
 &\leq
 \widehat{\mathbb E}
 \left[
  g_\delta(\widehat X_t-\widehat X_s)
 \right]  \\
 &=\widetilde{\mathbb E}
   \left[\Phi_h^\delta\right].
\end{split}
\]
Since
\[
        \widetilde{\mathbb E}[V]<\infty,
\]
it follows that
\[
 \widehat{\mathbb E}
 \left[
  I_{\{|\widehat X_t-\widehat X_s|\geq\delta\}}
 \right]
 \leq
 \frac{4C}{\delta^2}
 \left(1+\widetilde{\mathbb E}[V]\right)(t-s).
\]
The right-hand side tends to zero as \(t\) decreases to \(s\).
Therefore the canonical stationary process \(\widehat X\) is
stochastically continuous.
We next verify the fixed-space condition in Theorem 3.22 of
Feng--Zhao. Let \(\varphi\in B_b(\mathbb R^d)\) satisfy
\[
        T_t\varphi=\varphi,
        \qquad
        T_t(-\varphi)=-\varphi,
        \qquad t\geq0.
\]
The first identity and Proposition~\ref{prop:fixed} imply that
\[
        \varphi(x)
        =
        \lim_{t\to\infty}T_t\varphi(x)
        =
        \widetilde{\mathbb E}[\varphi],
        \qquad x\in\mathbb R^d.
\]
Thus \(\varphi\) is constant on \(\mathbb R^d\), and hence is
constant \(\widetilde{\mathbb E}\)-quasi-surely. Condition (i) of
Theorem 3.22 of Feng--Zhao is therefore satisfied.
Since \(\widehat X\) is stochastically continuous, the implication
\(\textup{(i)}\Rightarrow\textup{(ii)}\) in that theorem shows that
\[
 \bigl(
   \widehat\Omega,\widehat{\mathcal D},
   (\widehat\theta_t)_{t\in\mathbb R},
   \widehat{\mathbb E}_{\widetilde{\mathbb E}}
 \bigr)
\]
is ergodic in the sense of Definition~\ref{def:ergodicity}. Accordingly,
\(\widetilde{\mathbb E}\) is an ergodic invariant sublinear
expectation.
The additional no-mean-uncertainty hypothesis in Theorem 3.22 of
Feng--Zhao is needed only for the converse implication and is not
required here.
\end{proof}
\appendix
\section{Extension to the weighted space}
\label{sec:weighted-extension}
\begin{lemma}
\label{lem:BbBVbridge}
Consider either the multiplicative semigroup \eqref{eq:multT} under
Assumption~\ref{ass:multG}, or the full semigroup \eqref{eq:fullT} under
Assumption~\ref{ass:fullG}. Let
\[
 V(x)=\abs{x}^2.
\]
For either model, denote its canonical representing family by
$\mathcal P_{\Q}$. Then:
\begin{enumerate}[label=\textup{(\roman*)}]
 \item for every $p\geq2$, $T>0$, and compact set $K\subset\R^d$,
 \begin{equation}\label{eq:BVmoment}
  \sup_{x\in K}\sup_{P\in\mathcal P_{\Q}}
  E^P\!\left[\sup_{0\leq r\leq T}\abs{X_r^x}^p\right]
  <\infty;
 \end{equation}
 \item for every $\varphi\in B_V(\R^d)$, the formula
 \begin{equation}\label{eq:BVextension}
  T_t\varphi(x)
  :=\sup_{P\in\mathcal P_{\Q}}
  E^P\!\left[\varphi(X_t^x)\right]
 \end{equation}
 is finite and
 \begin{equation}\label{eq:BVgrowth}
  \abs{T_t\varphi(x)}
  \leq
  \|\varphi\|_V
  \left[
   1+e^{-ct}V(x)+\frac{C_0}{c}(1-e^{-ct})
  \right],
 \end{equation}
 where $c,C_0$ are the constants in \eqref{eq:multGLyap} or
 \eqref{eq:fullLyap}. In particular,
 $T_t(B_V(\R^d))\subset B_V(\R^d)$;
 \item if
 \[
  \varphi_n:=(-n)\vee(\varphi\wedge n),
 \]
 then, for every $T>0$ and compact $K\subset\R^d$,
 \begin{equation}\label{eq:BVtruncation}
  \lim_{n\to\infty}
  \sup_{\substack{0\leq t\leq T\\x\in K}}
  \abs{T_t\varphi_n(x)-T_t\varphi(x)}
  =0;
 \end{equation}
 \item the dynamic programming identity extends to $B_V(\R^d)$:
 \begin{equation}\label{eq:BVDPP}
  T_{t+s}\varphi=T_t(T_s\varphi),
  \qquad s,t\geq0,\quad\varphi\in B_V(\R^d);
 \end{equation}
 \item for every $\varphi\in B_V(\R^d)$, the function
 $u(t,x)=T_t\varphi(x)$ is locally bounded and locally H\"older continuous
 on $(0,\infty)\times\R^d$. It satisfies the same upper Bellman equation as
 in Proposition~\ref{prop:multbellman}, respectively
 Proposition~\ref{prop:fullHarnack}. For every constant $M$, the function
 $M-u$ satisfies the corresponding lower Bellman equation. Consequently,
 the local Harnack and oscillation arguments used in
 Theorems~\ref{thm:multG} and \ref{thm:fullG} remain valid for
 $\varphi\in B_V(\R^d)$.
\end{enumerate}
\end{lemma}
\begin{proof}
Under either set of assumptions, compactness of $\Q$ and the global
Lipschitz and linear-growth bounds imply that the drift and diffusion
coefficients are globally Lipschitz and have linear growth uniformly over
$P\in\mathcal P_{\Q}$. The standard BDG and Gronwall estimates therefore
give
\[
 E^P\!\left[\sup_{0\leq r\leq T}\abs{X_r^x}^p\right]
 \leq C_{p,T}(1+\abs{x}^p),
\]
where $C_{p,T}$ is independent of $P$. This proves
\eqref{eq:BVmoment}.
Let $\varphi\in B_V(\R^d)$. By monotonicity and sublinearity,
\[
 \abs{T_t\varphi(x)}
 \leq T_t\abs{\varphi}(x)
 \leq\|\varphi\|_V\bigl(1+T_tV(x)\bigr).
\]
It\^o's formula under each $P\in\mathcal P_{\Q}$ and the corresponding
robust Lyapunov condition give
\[
 T_tV(x)
 \leq e^{-ct}V(x)+\frac{C_0}{c}(1-e^{-ct}).
\]
This proves \eqref{eq:BVgrowth} and shows that $T_t\varphi\in B_V(\R^d)$.
We next prove \eqref{eq:BVtruncation}. The elementary inequality
$\abs{\sup_P a_P-\sup_P b_P}\leq\sup_P\abs{a_P-b_P}$ gives
\[
 \abs{T_t\varphi_n(x)-T_t\varphi(x)}
 \leq
 \sup_{P\in\mathcal P_{\Q}}
 E^P\!\left[
  \abs{\varphi_n(X_t^x)-\varphi(X_t^x)}
 \right].
\]
Moreover,
\[
 \abs{\varphi_n(y)-\varphi(y)}
 \leq
 \|\varphi\|_V(1+\abs{y}^2)
 I_{\{\abs{\varphi(y)}>n\}}.
\]
Since
$\abs{\varphi(y)}\leq\|\varphi\|_V(1+\abs{y}^2)$, there are
$R_n\uparrow\infty$ such that
\[
 \{\abs{\varphi(y)}>n\}\subset\{\abs{y}>R_n\}.
\]
Choose $\delta>0$. For $R\geq1$,
\[
 (1+\abs{y}^2)I_{\{\abs{y}>R\}}
 \leq 2R^{-\delta}\abs{y}^{2+\delta}.
\]
The moment estimate \eqref{eq:BVmoment}, with $p=2+\delta$, now yields
\[
 \sup_{\substack{0\leq t\leq T\\x\in K}}
 \sup_{P\in\mathcal P_{\Q}}
 E^P\!\left[
  \abs{\varphi_n(X_t^x)-\varphi(X_t^x)}
 \right]
 \longrightarrow0,
\]
which proves \eqref{eq:BVtruncation}.
For bounded $\varphi_n$, the bounded-data dynamic programming principle
gives
\begin{equation}\label{eq:BVDPPbounded}
 T_{t+s}\varphi_n=T_t(T_s\varphi_n).
\end{equation}
Set $g_n=T_s\varphi_n$ and $g=T_s\varphi$. By
\eqref{eq:BVtruncation}, $g_n\to g$ uniformly on compact sets, while
\eqref{eq:BVgrowth} gives
\[
 \abs{g_n(y)}+\abs{g(y)}
 \leq C_s\|\varphi\|_V(1+V(y))
\]
uniformly in $n$. Hence, for a compact $K\subset\R^d$ and $R\geq1$,
\[
 \begin{split}
 \sup_{x\in K}\abs{T_tg_n(x)-T_tg(x)}
 &\leq
 \sup_{\abs{y}\leq R}\abs{g_n(y)-g(y)}
 \\
 &\quad
 +C_s\|\varphi\|_V
 \sup_{\substack{x\in K\\P\in\mathcal P_{\Q}}}
 E^P\!\left[
  (1+\abs{X_t^x}^2)I_{\{\abs{X_t^x}>R\}}
 \right].
 \end{split}
\]
First let $n\to\infty$ and then $R\to\infty$, using
\eqref{eq:BVmoment}. Passing to the limit in
\eqref{eq:BVDPPbounded} proves \eqref{eq:BVDPP}.
Finally, let $u_n(t,x)=T_t\varphi_n(x)$. The bounded-data results in
Sections~\ref{sec:multG} and~\ref{sec:fullG} show that $u_n$ is locally H\"older continuous on
positive-time cylinders and satisfies the appropriate upper Bellman
equation. By \eqref{eq:BVtruncation}, $u_n\to u$ locally uniformly on
$[0,T]\times K$ for every $T>0$ and compact $K$. Stability of viscosity
solutions under locally uniform convergence therefore shows that $u$
satisfies the same Bellman equation on $(0,\infty)\times\R^d$.
Estimate \eqref{eq:BVgrowth} makes $u$ bounded on every compact
positive-time cylinder. The local interior oscillation estimate for the
common Pucci class then gives local H\"older continuity. Since the Bellman
operators annihilate constants, $M-u$ satisfies the corresponding lower
equation. The local Harnack and oscillation arguments therefore apply to
$T_t\varphi$ for every $\varphi\in B_V(\R^d)$.
\end{proof}

\bibliographystyle{amsplain}
\bibliography{refs}
\end{document}